# Analysis and Control of a Non-Standard Hyperbolic PDE Traffic Flow Model

**Iasson Karafyllis[*], Nikolaos Bekiaris-Liberis[**] and Markos Papageorgiou[**]**

[*]Dept. of Mathematics, National Technical University of Athens, Zografou Campus, 15780, Athens, Greece, email: iasonkar@central.ntua.gr

[**] Dynamic Systems and Simulation Laboratory, Technical University of Crete, Chania, 73100, Greece
(emails: nikos.bekiaris@dssl.tuc.gr, markos@dssl.tuc.gr)

**Abstract**

The paper provides results for a non-standard, hyperbolic, 1-D, nonlinear traffic flow model on a bounded domain. The model consists of two first-order PDEs with a dynamic boundary condition that involves the time derivative of the velocity. The proposed model has features that are important from a traffic-theoretic point of view: is completely anisotropic and information travels forward exactly at the same speed as traffic. It is shown that, for all physically meaningful initial conditions, the model admits a globally defined, unique, classical solution that remains positive and bounded for all times. Moreover, it is shown that global stabilization can be achieved for arbitrary equilibria by means of an explicit boundary feedback law. The stabilizing feedback law depends only on the inlet velocity and consequently, the measurement requirements for the implementation of the proposed boundary feedback law are minimal. The efficiency of the proposed boundary feedback law is demonstrated by means of a numerical example.

**Keywords:** hyperbolic PDEs, traffic flow, boundary feedback.

## 1. Introduction

The study of vehicular traffic flow by means of hyperbolic Partial Differential Equations (PDEs) started in the 1950s with the LWR first-order model (see [26,31]). In order to describe more accurately the velocity dynamics, second-order models were later studied (see [1,28,36]). All 1-D traffic flow models were developed for unbounded domains (usually the whole real axis). Researchers working on second-order models as well as critics of second-order models (see [11]) have agreed that a valid traffic flow model must: (i) include the vehicle conservation equation, (ii) admit bounded solutions which predict positive values for both density and velocity, (iii) obey the so-called anisotropy principle, i.e., the fact that a vehicle is influenced only by the traffic dynamics ahead of it, (iv) not allow waves traveling forward with a speed greater than the traffic speed. Recently, researchers have developed two phase models (see [7,24]), which agree with experimental results that report strong differences between the free and congested vehicular flow.

Recent advances in the boundary feedback control of hyperbolic systems of PDEs (see for instance [2,3,6,8,9,10,12,13,18,21,22,29,30,34,35]) as well as advances in the control of discrete-time, finite-dimensional traffic flow models (see [16,17,19,27] and references therein) have motivated the study of well-posedness and control of traffic flow models on bounded domains. Both issues (well-posedness and control) for first-order models in bounded domains were studied in [4,5,32] by means of boundary conditions at the inlet and outlet that may or may not become active at certain time instants. The stabilization of equilibrium profiles for second-order models in bounded domains by means of boundary feedback was also studied in [23,37].

The present work presents a novel, hyperbolic, nonlinear, second-order, 1-D traffic flow model on a bounded domain. The arguments leading to the derivation of the model are based on the assumption that the road is relatively crowded. It consists of two quasilinear first-order PDEs with a dynamic nonlinear boundary condition that involves the time derivative of the velocity, which may be viewed as boundary relaxation, analogously to in-domain relaxation in second-order traffic flow models [1,36]. The presence of this dynamic boundary condition makes the model non-standard, and thus, the existence and uniqueness of its solutions cannot be guaranteed by using standard results (see [2,20,25]). The existence and uniqueness issues are studied in the present work. Specifically, it is shown that for all physically meaningful initial conditions the model admits a globally defined, unique, classical solution that remains positive and bounded for all times. As a result, we can guarantee that the proposed model has all of the four features mentioned in the first paragraph that are important from a traffic-theoretic point of view. The second contribution of the present work is the study of the control problem for the proposed model. Specifically, we design a simple, nonlinear, boundary feedback law, adjusting the inlet flow (via, e.g., ramp metering). The boundary feedback law employs only measurements of the inlet velocity, and consequently, the measurement requirements for implementation of the proposed controller are minimal. Moreover, it is shown that the developed control design achieves global asymptotic stabilization of arbitrary equilibria, in the sup-norm of the logarithmic deviation of the state from its equilibrium point. The efficiency of the proposed feedback law is demonstrated by means of a numerical example.

The structure of the present work is as follows: Section 2 is devoted to the presentation of the model and the statement of the first main result (Theorem 2.1) which guarantees, for all physically meaningful initial conditions, the existence of a globally defined, unique, classical solution that remains positive and bounded for all times. The control design and the statement of the second main result, which guarantees global stabilization of arbitrary equilibria of the model (Theorem 3.1) are given in Section 3. A simple illustrating example is presented in Section 4. The proofs of the main results as well as the statement of two auxiliary results are provided in Section 5. One of the auxiliary results has interest on its own (Proposition 5.2), because it covers a case not studied in [2,20,25]: the case of a transport PDE with a non-negative (possibly zero at some points) transport velocity. A unique, classical solution is shown to exist, which is differentiable and satisfies the PDE even on the boundary (something that cannot be guaranteed by the results in [20]). The concluding remarks are provided in Section 6. Finally, the Appendix contains the proofs of the two auxiliary results, which were stated in Section 5.

**Notation.** Throughout this paper, we adopt the following notation.

* $\Re_+ := [0,+\infty)$. For a real number $x \in \Re$, $[x]$ denotes the integer part of $x$, i.e., the greatest integer which is less or equal to $x$.
* Let $U \subseteq \Re^n$ be a set with non-empty interior and let $\Omega \subseteq \Re$ be a set. By $C^0(U;\Omega)$, we denote the class of continuous mappings on $U$, which take values in $\Omega$. By $C^k(U;\Omega)$, where $k \geq 1$, we denote the class of continuous functions on $U$, which have continuous derivatives of order $k$ on $U$ and take values in $\Omega$. When $\Omega$ is omitted, i.e., when we write $C^k(U)$, it is implied that $\Omega = \Re$.
* Let $T \in (0,+\infty)$ and $u:[0,T]\times[0,1]\to\Re$ be given. We use the notation $u[t]$ to denote the profile at certain $t \in [0,T]$, i.e., $(u[t])(x) = u(t,x)$ for all $x \in [0,1]$. For a bounded $w:[0,1]\to\Re$ the sup-norm is given by $\|w\|_\infty := \sup_{0 \leq x \leq 1}(|w(x)|)$.
* $W^{2,\infty}([0,1])$ is the Sobolev space of $C^1$ functions on $[0,1]$ with Lipschitz derivative.
* By $K$ we denote the class of strictly increasing continuous functions $a:\Re_+ \to \Re_+$ with $a(0) = 0$. By $K_\infty$ we denote the class of functions $a \in K$ with $\lim_{s \to +\infty} a(s) = +\infty$. By $KL$ we denote the set of

all functions $\sigma \in C^0(\Re_+ \times \Re_+ ; \Re_+)$ with the properties: (i) for each $t \geq 0$, $\sigma(\cdot, t)$ is of class $K$; (ii) for each $s \geq 0$, $\sigma(s, \cdot)$ is non-increasing with $\lim_{t \to +\infty} \sigma(s,t) = 0$.

## 2. A Non-Standard Traffic Flow Model

*2.I. Model Description*

Second-order traffic flow models involve a system of hyperbolic PDEs on the positive semiaxis. The state variables are the vehicle density $\rho(t,x)$ and the vehicle velocity $v(t,x)$, where $t \geq 0$ is time and $x$ is the spatial variable. All traffic flow models involve the conservation equation

$$\frac{\partial \rho}{\partial t}(t,x) + v(t,x)\frac{\partial \rho}{\partial x}(t,x) + \rho(t,x)\frac{\partial v}{\partial x}(t,x) = 0 \tag{2.1}$$

and an additional PDE for the velocity. In a relatively crowded road, the vehicle velocity depends heavily on the velocity of downstream vehicles. Therefore, the following equation may be appropriate for the description of the evolution of the velocity profile:

$$\frac{\partial v}{\partial t}(t,x) - c\frac{\partial v}{\partial x}(t,x) = 0, \tag{2.2}$$

where $c > 0$ is a constant related to the drivers' speed of adjusting their velocity. Equation (2.2) may also arise as a linearization of the equation of the Aw-Rascle-Zhang model (see [1,36]) without an in-domain relaxation term. Here, we consider the model (2.1), (2.2) on a bounded domain, i.e., we assume that $x \in [0,1]$. The full model requires the specification of two boundary conditions. One boundary condition describes the inlet conditions and more particularly the effect of the inlet demand $q(t) > 0$ and takes the form

$$\rho(t,0) = h\left(\frac{q(t)}{v(t,0)}\right), \text{ for } t \geq 0 \tag{2.3}$$

where $h \in C^2(\Re_+)$ is a non-decreasing function that satisfies

$$h(s) = s \text{ for } s \in [0, \rho_{\max} - \varepsilon] \text{ and } h(s) = \rho_{\max} \text{ for } s \geq \rho_{\max}, \tag{2.4}$$

where $\rho_{\max} > 0$ is a constant related to the physical upper bound of density in the particular road and $\varepsilon \in (0, \rho_{\max})$ is a sufficiently small constant. Notice that (2.3) implies that the inlet demand $q(t) > 0$ is equal to the vehicle inflow $\rho(t,0)v(t,0)$, provided that $q(t) \leq (\rho_{\max} - \varepsilon)v(t,0)$. The boundary condition (2.3) as well as the rest of the model (2.1), (2.2) comes together with the following requirement:
$$\rho(t,x) > 0 \text{ and } v(t,x) > 0, \text{ for all } (t,x) \in \Re_+ \times [0,1] \tag{2.5}$$

Condition (2.5) is an essential requirement for traffic flow models and it should be noticed here that some second-order traffic flow models do not meet this requirement. In what follows, we show that the proposed model meets this requirement.

In order to have a well-posed hyperbolic system, we also need a boundary condition at the outlet $x=1$. Assuming that the flow downstream the outlet is uncongested (free), it is reasonable to assume that the relaxation term becomes dominant. So, we get

$$\frac{\partial v}{\partial t}(t,1) = -\mu\big(v(t,1) - f(\rho(t,1))\big), \text{ for } t \geq 0 \tag{2.6}$$

where $\mu > 0$ is a constant and $f \in C^1(\Re_+)$ is a positive, bounded, non-increasing function that expresses the fundamental diagram relation between density and velocity.

*2.II. Traffic-Theoretic Features of the Model*

Equations (2.1), (2.2), (2.3), (2.6) form a non-standard system of nonlinear hyperbolic PDEs. The reason that system (2.1), (2.2), (2.3), (2.6) cannot be studied by existing results in hyperbolic systems (see [2,20,25]) is the non-standard boundary condition (2.6). However, in what follows, we show that system (2.1), (2.2), (2.3), (2.6) exhibits unique, positive, globally defined $C^1$ solutions for all positive initial conditions. Moreover, we show that density and velocity are bounded from above by certain bounds that depend only on the initial conditions and the physical upper bounds of the density and velocity, i.e., $\rho_{max}$ and $v_{max} = f(0)$, respectively. Before we show this, it is important to emphasize that (2.1), (2.2), (2.3), (2.6):

- is a traffic flow model that can be applied to bounded domains, i.e., $x \in [0,1]$, without imposing a boundary condition with no physical meaning or assuming knowledge of the density/velocity out of the domain,

- is completely anisotropic, i.e., the velocity depends only on the velocity of downstream vehicles,

- is a hyperbolic model with two eigenvalues $v$ and $-c$; consequently, information travels forward exactly in the same speed as traffic,

- allows only equilibria which satisfy the fundamental diagram law $v = f(\rho)$, i.e., when $q(t) \equiv q_{eq} > 0$ then the equilibrium profiles are given by $\rho(x) \equiv \rho_{eq}$, $v(x) \equiv f(\rho_{eq})$, where $\rho_{eq} > 0$ is a solution of $\rho_{eq} = h\left(\dfrac{q_{eq}}{f(\rho_{eq})}\right)$.

All the above features are important for a traffic flow model.

*2.III. Characteristic Form of the System*

Let $\rho_{eq} \in (0, \rho_{max})$ be a given constant. The nonlinear transformation of the density variable

$$\rho(t,x) = \rho_{eq} \exp(w(t,x)) \frac{c + f(\rho_{eq})}{c + v(t,x)} \qquad (2.7)$$

gives the equation

$$\frac{\partial w}{\partial t}(t,x) + v(t,x) \frac{\partial w}{\partial x}(t,x) = 0 \qquad (2.8)$$

with the boundary conditions

$$w(t,0) = \ln\left(\rho_{eq}^{-1} h\left(\frac{q(t)}{v(t,0)}\right) \frac{c + v(t,0)}{c + f(\rho_{eq})}\right), \quad \frac{\partial v}{\partial t}(t,1) = -\mu\left(v(t,1) - f\left(\rho_{eq} \exp(w(t,1)) \frac{c + f(\rho_{eq})}{c + v(t,1)}\right)\right) \qquad (2.9)$$

The hyperbolic system (2.2), (2.8), (2.9) is nothing else but the hyperbolic system (2.1), (2.2), (2.3), (2.6) in Riemann coordinates. Provided that the initial conditions are positive, i.e., $\rho(0,x) > 0$, $v(0,x) > 0$, for $x \in [0,1]$, we are in a position to construct a unique solution to (2.1), (2.2), (2.3), (2.6) by constructing a unique solution to (2.2), (2.8), (2.9) and employing the nonlinear transformation (2.7).

*2.IV. First Main Result*

The solution of (2.2), (2.8), (2.9) is constructed by the following theorem. Its proof is provided in Section 5.

**Theorem 2.1:** *Let $a \in C^2(\Re_+ \times \Re_+)$ be a given function and let $c > 0$, $\mu \geq 0$ be given constants. Let $g \in C^1(\Re_+ \times \Re)$ be a given function for which there exists a constant $v_{\max} > 0$ such that the following inequality holds*

$$0 < g(0, w) \leq g(v, w) \leq v_{\max}, \text{ for all } v \in \Re_+, w \in \Re \quad (2.10)$$

*Let $\theta, \varphi \in W^{2,\infty}([0,1])$ be given functions with $\varphi(x) > 0$ for all $x \in [0,1]$, for which the equalities $\theta(0) = a(0, \varphi(0))$, $\frac{\partial a}{\partial t}(0, \varphi(0)) + c \frac{\partial a}{\partial v}(0, \varphi(0)) \varphi'(0) = -\varphi(0) \theta'(0)$, $\varphi'(1) = -\mu c^{-1} (\varphi(1) - g(\varphi(1), \theta(1)))$ hold. Then the initial-boundary value problem*

$$\frac{\partial w}{\partial t}(t, x) + v(t, x) \frac{\partial w}{\partial x}(t, x) = \frac{\partial v}{\partial t}(t, x) - c \frac{\partial v}{\partial x}(t, x) = 0, \text{ for all } (t, x) \in \Re_+ \times [0,1] \quad (2.11)$$

$$w(t, 0) - a(t, v(t, 0)) = \frac{\partial v}{\partial t}(t, 1) + \mu(v(t, 1) - g(v(t, 1), w(t, 1))) = 0, \text{ for all } t \geq 0 \quad (2.12)$$

$$w(0, x) - \theta(x) = v(0, x) - \varphi(x) = 0, \text{ for all } x \in [0,1] \quad (2.13)$$

*admits a unique solution $w, v \in C^1(\Re_+ \times [0,1])$. Moreover, the solution $w, v \in C^1(\Re_+ \times [0,1])$ has Lipschitz derivatives on every compact $S \subset \Re_+ \times [0,1]$ and satisfies the following inequalities for all $(t, x) \in \Re_+ \times [0,1]$:*

$$\|w[t]\|_\infty \leq \max(B_t, \|\theta\|_\infty) \quad (2.14)$$

$$\min\left(\min_{0 \leq x \leq 1}(\varphi(x)), \min\left\{g(0, w) : |w| \leq \max(B_t, \|\theta\|_\infty)\right\}\right) \leq v(t, x) \leq \max\left(\max_{0 \leq x \leq 1}(\varphi(x)), v_{\max}\right) \quad (2.15)$$

*where $B_t := \max\left\{|a(s, v)| : s \in [0, t], 0 \leq v \leq \max\left(\max_{0 \leq x \leq 1}(\varphi(x)), v_{\max}\right)\right\}$.*

**Remark 2.2:** Theorem 2.1 shows that the appropriate space (state space) for studying the hyperbolic system (2.1), (2.2), (2.3), (2.6) is the space

$$X = \left\{ (\rho, v) \in \left(W^{2,\infty}([0,1])\right)^2 : \begin{array}{c} \min(\rho(x), v(x)) > 0 \text{ for all } x \in [0,1], \\ cv'(1) = -\mu(v(1) - f(\rho(1))), \\ \exists a_1 > 0, a_2 \in \Re \text{ such that } \rho(0) = h(a_1), v(0)\rho'(0) + \rho(0)v'(0) = a_2 h'(a_1) \end{array} \right\}. \quad (2.16)$$

In order to construct a solution $(\rho[t], v[t]) \in X$ of (2.1), (2.2), (2.3), (2.6) with initial conditions in $(\rho_0, v_0) \in X$, we apply Theorem 2.1 with

$$a(t, v) := \begin{cases} \ln\left(\rho_{eq}^{-1} h\left(\frac{q(t)}{v}\right) \frac{c + v}{c + f(\rho_{eq})}\right) & \text{if } v > 0 \\ \ln\left(\rho_{eq}^{-1} \rho_{\max} \frac{c + v}{c + f(\rho_{eq})}\right) & \text{if } v = 0 \end{cases}, \quad g(v, w) := f\left(\rho_{eq} \exp(w) \frac{c + f(\rho_{eq})}{c + v}\right), \quad v_{\max} := f(0),$$

$$\theta(x) = \ln\left(\frac{\rho_0(x)(c + v_0(x))}{(c + f(\rho_{eq}))\rho_{eq}}\right), \quad \varphi(x) = v_0(x) \text{ for all } x \in [0,1]$$

and we consider $q \in C^2(\Re_+; (0, +\infty))$ to be the input of the model. The set of admissible inputs consists of all functions $q \in C^2(\Re_+; (0, +\infty))$ that satisfy the compatibility conditions $\rho_0(0) = h\left(\frac{q(0)}{v_0(0)}\right)$ and $v_0(0)\rho_0'(0) + \rho_0(0)v_0'(0) + h'\left(\frac{q(0)}{v_0(0)}\right)\frac{\dot{q}(0)}{v_0(0)} = ch'\left(\frac{q(0)}{v_0(0)}\right)\frac{q(0)}{v_0^2(0)}v_0'(0)$. The solution $(\rho[t], v[t]) \in X$ of (2.1), (2.2), (2.3), (2.6) is found by using the solution $(w[t], v[t])$ of (2.11), (2.12), (2.13) in conjunction with formula (2.7). Notice that if $v_0(x) \leq v_{\max}$ for all $x \in [0,1]$, then estimate (2.15) implies that $0 < v(t,x) \leq v_{\max}$ for all $(t,x) \in \Re_+ \times [0,1]$ and for all admissible $q \in C^2(\Re_+;(0,+\infty))$. Similarly, by performing more detailed calculations than those in the proof of Theorem 2.1, we are in a position to verify that if $\rho_0(x) \leq \rho_{\max}\frac{c+v_{\max}}{c}$ for all $x \in [0,1]$, then the estimate $0 < \rho(t,x) \leq \rho_{\max}\frac{c+v_{\max}}{c}$ holds for all $(t,x) \in \Re_+ \times [0,1]$ and for all admissible $q \in C^2(\Re_+;(0,+\infty))$. When $c \gg v_{\max}$, the previous estimate implies that the upper bound for density is approximately $\rho_{\max}$.

## 3. Controlling the Traffic Flow Model

*3.I. Motivation for Control Design*

The fact that for the case $q(t) \equiv q_{eq} > 0$ the equilibrium profiles for (2.1), (2.2), (2.3), (2.6) are given by $\rho(x) \equiv \rho_{eq}$, $v(x) \equiv f(\rho_{eq})$, where $\rho_{eq} > 0$ is a solution of $\rho_{eq} = h\left(\frac{q_{eq}}{f(\rho_{eq})}\right)$, implies that there may be multiple equilibria. For example, for the case $f(\rho) = A\exp(-b\rho)$, where $A, b > 0$ are constants (that corresponds to the so-called Underwood model; see for instance [33]) if $q_{eq} \in \left[A\rho_{\max}\exp(-b\rho_{\max}), \frac{A}{b}\exp(-1)\right]$ and $1 < b\rho_{\max}$ then there are (at least) two solutions of the equation $\rho_{eq} = h\left(\frac{q_{eq}}{f(\rho_{eq})}\right)$: one solution in the interval $(0, b^{-1}]$ and $\rho_{\max}$. Consequently, it is not possible to guarantee that for every initial condition $(\rho_0, v_0) \in X$ with $\rho_0(0) = h\left(\frac{q_{eq}}{v_0(0)}\right)$, $v_0(0)\rho_0'(0) + \rho_0(0)v_0'(0) = ch'\left(\frac{q(0)}{v_0(0)}\right)\frac{q_{eq}}{v_0^2(0)}v_0'(0)$, the solution $(\rho[t], v[t]) \in X$ of (2.1), (2.2), (2.3), (2.6) with $q(t) \equiv q_{eq} > 0$ will converge to a specific equilibrium profile as $t \to +\infty$. This implies lack of global asymptotic stability. Moreover, such cases are the ones that ideally one would like to have: for the case $f(\rho) = A\exp(-b\rho)$, where $A, b > 0$ are constants, the ideal operation of the freeway would be exactly where the flow becomes maximized, i.e., when $\rho = b^{-1}$. Notice that in this case and if $1 \leq b(\rho_{\max} - \varepsilon)$, where $\varepsilon \in (0, \rho_{\max})$ is the constant involved in (2.4), $q_{eq} = \frac{A}{b}\exp(-1)$ and we have (at least) two equilibria: $b^{-1}$ and $\rho_{\max}$. In such cases, global stabilization of a specific equilibrium profile may be achieved by boundary feedback control.

*3.II. Collocated Boundary Control Design and Stability Analysis*

The following theorem shows that stabilization of the equilibrium profile for a given desired equilibrium density $\rho_{eq} > 0$ can be achieved by controlling the inlet flow. It is important to notice that the stabilizing feedback law depends *only* on the inlet velocity. Therefore, the measurement requirements for the implementation of the proposed boundary feedback law are minimal.

**Theorem 3.1:** *Consider the nonlinear traffic flow model (2.1), (2.2), (2.3), (2.6) and let $\rho_{eq} > 0$ be the desired equilibrium density. Suppose that $\rho_{eq} \leq \frac{c}{c+f(\rho_{eq})}(\rho_{\max} - \varepsilon)$ and that the following inequality holds:*

$$\left(v - f\left(\rho_{eq} \frac{c+f(\rho_{eq})}{c+v}\right)\right)\left(v - f(\rho_{eq})\right) > 0, \text{ for all } v \geq 0, v \neq f(\rho_{eq}) \quad (3.1)$$

*Then there exists a function $Q \in KL$ such that for every $(\rho_0, v_0) \in X$ for which the equalities $\rho_0(0) = \rho_{eq} \frac{c+f(\rho_{eq})}{c+v_0(0)}$, $\rho_0'(0) = -\frac{\rho_0(0)}{c+v_0(0)} v_0'(0)$ hold, the initial-boundary value problem (2.1), (2.2), (2.3), (2.6) with*

$$q(t) = \rho_{eq} v(t,0) \frac{c+f(\rho_{eq})}{c+v(t,0)} \quad (3.2)$$

$$\rho(0,x) - \rho_0(x) = v(0,x) - v_0(x) = 0, \text{ for all } x \in [0,1] \quad (3.3)$$

*admits a unique solution $\rho, v \in C^1(\Re_+ \times [0,1])$, with $(\rho[t], v[t]) \in X$ for all $t \geq 0$ satisfying the following estimate for all $t \geq 0$:*

$$\max_{0 \leq x \leq 1}\left(\left|\ln\left(\frac{\rho(t,x)}{\rho_{eq}}\right)\right|\right) + \max_{0 \leq x \leq 1}\left(\left|\ln\left(\frac{v(t,x)}{f(\rho_{eq})}\right)\right|\right) \leq Q\left(\max_{0 \leq x \leq 1}\left(\left|\ln\left(\frac{\rho_0(x)}{\rho_{eq}}\right)\right|\right) + \max_{0 \leq x \leq 1}\left(\left|\ln\left(\frac{v_0(x)}{f(\rho_{eq})}\right)\right|\right), t\right) \quad (3.4)$$

**Remark 3.2:** Notice that inequality (3.1) holds automatically for $v \geq f(0)$ and $v = 0$. Thus, inequality (3.1) is equivalent to the following implications

$$\rho_{eq} > \rho > \rho_{eq} \frac{c+f(\rho_{eq})}{c+f(0)} \Rightarrow \rho_{eq}(c+f(\rho_{eq})) > \rho(c+f(\rho))$$

$$\rho_{eq} < \rho < \rho_{eq} \frac{c+f(\rho_{eq})}{c} \Rightarrow \rho_{eq}(c+f(\rho_{eq})) < \rho(c+f(\rho))$$

Therefore, a sufficient condition for (3.1) is the assumption that the function $F(\rho) := \rho(c+f(\rho))$ is increasing on the interval $\left(\rho_{eq} \frac{c+f(\rho_{eq})}{c+f(0)}, \rho_{eq} \frac{c+f(\rho_{eq})}{c}\right)$. Consequently, (3.1) holds automatically when $c + f(\rho) + \rho f'(\rho) > 0$ for all $\rho \in \left(\rho_{eq} \frac{c+f(\rho_{eq})}{c+f(0)}, \rho_{eq} \frac{c+f(\rho_{eq})}{c}\right)$. For example, when $f(\rho) = A\exp(-b\rho)$, where $A, b > 0$ are constants (Underwood model), we guarantee that (3.1) holds when the inequality $c\exp(b\rho) + A > Ab\rho$ holds for $\rho \in \left(\rho_{eq} \frac{c+A\exp(-b\rho_{eq})}{c+A}, \rho_{eq} \frac{c+A\exp(-b\rho_{eq})}{c}\right)$. It should be noticed that in this case (3.1) holds automatically when the velocity ratio $A/c$ is sufficiently small no matter what $\rho_{eq}$ is: when $c\exp(2) \geq A$ the function $F(\rho) := \rho(c+A\exp(-b\rho))$ is increasing on $\Re_+$.

**Remark 3.3:** When the compatibility conditions $\rho_0(0) = \rho_{eq}\frac{c+f(\rho_{eq})}{c+v_0(0)}$, $\rho_0'(0) = -\frac{\rho_0(0)}{c+v_0(0)} v_0'(0)$ do not hold, then we satisfy the compatibility conditions implied by Theorem 2.1, namely we find $a > 0$ and $b \in \Re$ so that $\rho_0(0) = h\left(\frac{a}{v_0(0)}\right)$ and $v_0(0)\rho_0'(0) + \rho_0(0)v_0'(0) = \frac{acv_0'(0) - bv_0(0)}{v_0^2(0)} h'\left(\frac{a}{v_0(0)}\right)$. In such a case, we must modify the control input so that the compatibility conditions hold; the control input can be given by the formula

$$q(t) = (1 - g_T(t))(a + bt) + g_T(t)\rho_{eq} v(t,0) \frac{c + f(\rho_{eq})}{c + v(t,0)}$$

where $T > 0$ is a small constant that satisfies $a + bT > 0$ and $g_T : \Re \to \Re$ is defined by $g_T(t) = 0$ for $t \leq 0$, $g_T(t) = 1$ for $t \geq T$ and $g_T(t) = \frac{\exp(-t^{-1})}{\exp(-t^{-1}) + \exp(-(T-t)^{-1})}$ for $t \in (0, T)$.

**Remark 3.4:** Estimate (3.4) is a stability estimate in the sup-norm of the logarithmic deviation of the state from its equilibrium values. The use of logarithmic deviation variables is customary for systems with positive state values (e.g., biological systems, see [18]).

## 4. Illustrative Example

We consider model (2.1), (2.2), (2.3), (2.6) with $f(\rho) = \frac{2}{5} \exp(1 - \rho)$ (Underwood model), $c = 5$, $\mu = 10$, $\rho_{\max} = 27/10$, $\varepsilon = 10^{-6}$, $h(s) = s(1 - g(s)) + \rho_{\max} g(s)$ for $s \geq 0$, where

$g(s) = 0$, for $s \in [0, \rho_{\max} - \varepsilon]$, $g(s) = 1$, for $s \geq \rho_{\max}$ and

$$g(s) = \frac{\exp\left(-(s + \varepsilon - \rho_{\max})^{-1}\right)}{\exp\left(-(s + \varepsilon - \rho_{\max})^{-1}\right) + \exp\left(-(\rho_{\max} - s)^{-1}\right)}, \text{ for } s \in (\rho_{\max} - \varepsilon, \rho_{\max}).$$

The objective is to stabilize the equilibrium point that maximizes the vehicle flow $\rho(x) \equiv \rho_{eq} = 1$, $v(x) \equiv f(\rho_{eq}) = 2/5$. It should be noticed that the open-loop system (2.1), (2.2), (2.3), (2.6) with $q(t) \equiv q_{eq} = 2/5$ has two equilibria: one is the desired equilibrium, and the other one is the fully congested equilibrium $\rho(x) \equiv \rho_{\max} = \frac{27}{10}$, $v(x) \equiv f(\rho_{\max}) = \frac{2}{5} \exp\left(-\frac{17}{10}\right)$. Numerical experiments show that the fully congested equilibrium attracts the solution of the open-loop system (2.1), (2.2), (2.3), (2.6) with $q(t) \equiv q_{eq} = \frac{2}{5}$ for many initial conditions. We choose the initial conditions

$\rho_0(x) = 1$ for $x \in [0, 9/20]$, $\rho_0(x) = 2$, for $x \in [1/2, 1]$,

$$\rho_0(x) = 1 + \frac{\exp\left(-(x - 9/20)^{-1}\right)}{\exp\left(-(x - 9/20)^{-1}\right) + \exp\left((x - 1/2)^{-1}\right)}, \text{ for } x \in \left(\frac{9}{20}, \frac{1}{2}\right), \text{ and}$$

$v_0(x) = f(\rho_0(x))$, for $x \in [0,1]$.

For this particular initial condition (but also for many others) the solution of the open-loop system (2.1), (2.2), (2.3), (2.6) with $q(t) \equiv q_{eq} = 2/5$ converges to the fully congested equilibrium $\rho(x) \equiv \rho_{\max} = \frac{27}{10}$, $v(x) \equiv f(\rho_{\max}) = \frac{2}{5} \exp\left(-\frac{17}{10}\right)$. The deviation of the solution from the desired equilibrium is shown in Fig. 1, where the evolution of the sup-norm of the logarithmic deviation from the desired equilibrium $X(t) := \max_{0 \leq x \leq 1}\left(\left|\ln\left(\frac{\rho(t,x)}{\rho_{eq}}\right)\right|\right) + \max_{0 \leq x \leq 1}\left(\left|\ln\left(\frac{v(t,x)}{f(\rho_{eq})}\right)\right|\right)$ is shown for the open-loop system (2.1), (2.2), (2.3), (2.6) with $q(t) \equiv q_{eq} = 2/5$.

In this case we can apply Theorem 3.1, since the condition $\rho_{eq} \leq \frac{c}{c + f(\rho_{eq})}(\rho_{\max} - \varepsilon)$ as well as condition (3.1) hold (recall Remark 3.2). Fig. 1 and Fig. 2 show the evolution of the density profile for the closed-loop system (2.1), (2.2), (2.3), (2.6) with (3.2). Fig. 3 and Fig. 4 show the convergence of the solution to the equilibrium profile $\rho(x) \equiv \rho_{eq} = 1$. It should be noted that at time $t = 6.58$, the solution has become identical (up to numerical accuracy) to the desired equilibrium. This is clear from Fig. 2, where it is shown the evolution of the sup-norm of the logarithmic deviation from the desired equilibrium $X(t) := \max_{0 \leq x \leq 1}\left(\left|\ln\left(\frac{\rho(t,x)}{\rho_{eq}}\right)\right|\right) + \max_{0 \leq x \leq 1}\left(\left|\ln\left(\frac{v(t,x)}{f(\rho_{eq})}\right)\right|\right)$ for the closed-loop system (2.1), (2.2), (2.3), (2.6) with (3.2). Fig. 5 shows the time evolution of the control input $q(t)$. The control input tries to keep the inlet density close to 1, while the heavy congestion belt is "washed out" slowly (due to small vehicle velocity in the congestion belt).

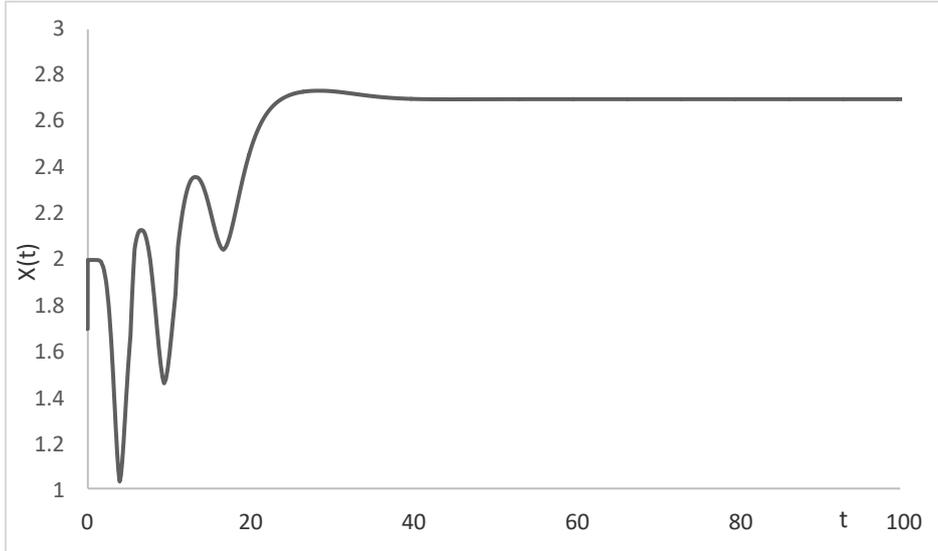

**Fig. 1:** Evolution of the sup-norm of the logarithmic deviation from the desired equilibrium $X(t) := \max_{0 \leq x \leq 1}\left(\left|\ln\left(\frac{\rho(t,x)}{\rho_{eq}}\right)\right|\right) + \max_{0 \leq x \leq 1}\left(\left|\ln\left(\frac{v(t,x)}{f(\rho_{eq})}\right)\right|\right)$ for the open-loop system (2.1), (2.2), (2.3), (2.6) with $q(t) \equiv q_{eq} = 2/5$.

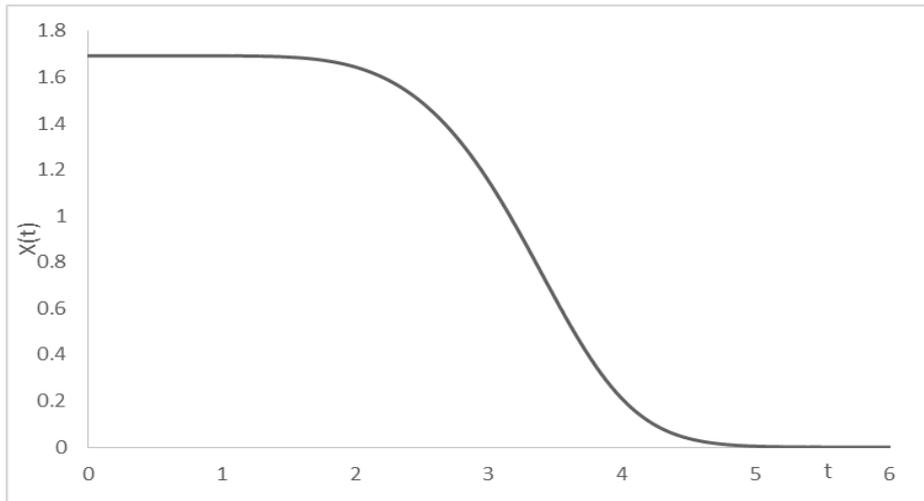

**Fig. 2:** Evolution of the sup-norm of the logarithmic deviation from the desired equilibrium $X(t) := \max_{0 \leq x \leq 1}\left(\left|\ln\left(\frac{\rho(t,x)}{\rho_{eq}}\right)\right|\right) + \max_{0 \leq x \leq 1}\left(\left|\ln\left(\frac{v(t,x)}{f(\rho_{eq})}\right)\right|\right)$ for the closed-loop system (2.1), (2.2), (2.3), (2.6) with (3.2).

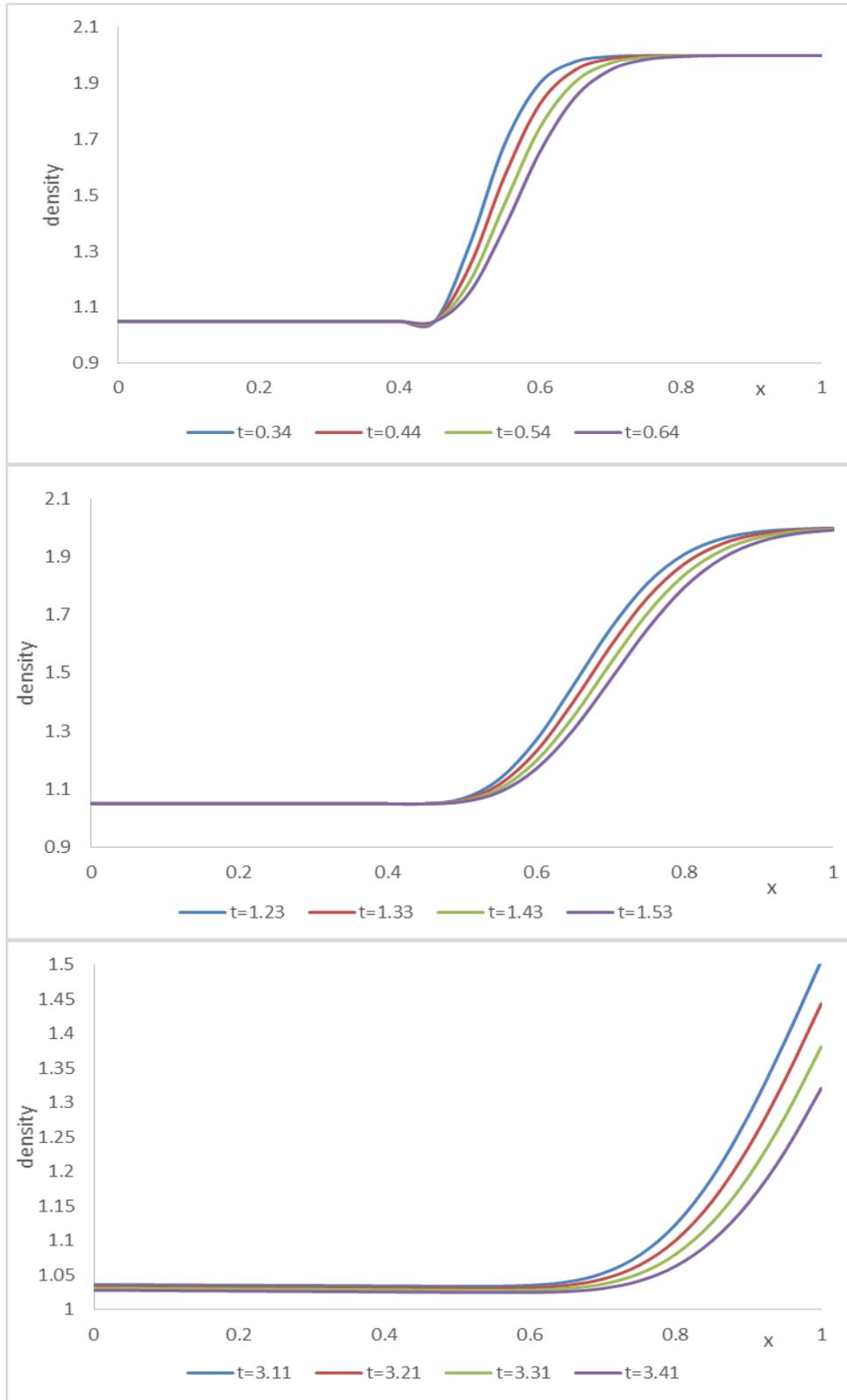

**Fig. 3:** Density profiles at various time instants for the closed-loop system (2.1), (2.2), (2.3), (2.6) with (3.2).

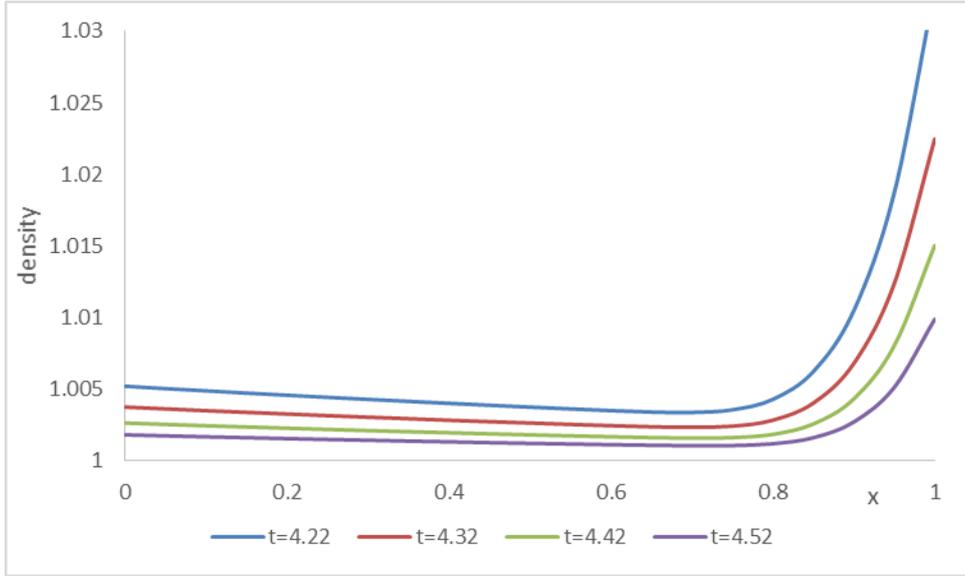

**Fig. 4:** Density profiles at various time instants for the closed-loop system (2.1), (2.2), (2.3), (2.6) with (3.2).

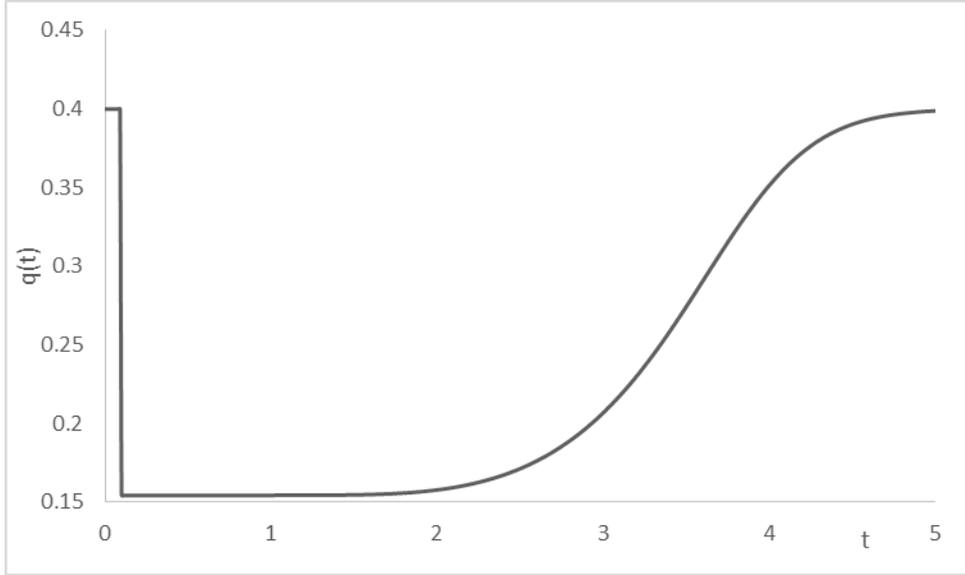

**Fig. 5:** The evolution of the control input $q(t)$ for the closed-loop system (2.1), (2.2), (2.3), (2.6) with (3.2).

## 5. Proofs of Main Results

*5.I. Technical Results*

The proof of Theorem 2.1 requires two technical results. Their proofs are given in the Appendix.

**Lemma 5.1:** *Suppose that there exist constants $a, b, p \geq 0$, $c > 0$ such that the sequence $\{x(k) \geq 0\}_{k=0}^{\infty}$ satisfies the inequality:*

$$x(k+1) \leq \max\left((1+a)x(k) + b, (1-c)x(k) + p\right), \text{ for all } k = 0,1,\ldots,m-1 \qquad (5.1)$$

*Then the following estimate holds:*

$$x(k) \leq \exp(ka)\left(x(0) + \frac{p}{a+c} + bk\right), \text{ for all } k = 0,1,\ldots,m \qquad (5.2)$$

The following auxiliary result has interest on its own, because it covers a case not studied in [2,20,25]: the case of a transport PDE with a non-negative (possibly zero at some points) transport velocity. A unique, classical solution is shown to exist, which is differentiable and satisfies the PDE even on the boundary (something that cannot be guaranteed by the results in [20]): this is important for the proof of Theorem 2.1, because uniform (Lipschitz) continuity of the derivatives of the solution on every compact set of the form $[0,T]\times[0,1]$ is used in an instrumental way.

**Proposition 5.2:** *Consider the initial-boundary value problem*

$$\frac{\partial w}{\partial t}(t,x) + v(t,x)\frac{\partial w}{\partial x}(t,x) = 0, \text{ for } t \geq 0,\ x \in [0,1] \tag{5.3}$$

$$w(0,x) = \varphi(x), \text{ for } x \in [0,1] \tag{5.4}$$

$$w(t,0) = a(t), \text{ for } t \geq 0 \tag{5.5}$$

*where* $\varphi \in W^{2,\infty}([0,1])$, $a \in W^{2,\infty}([0,T])$ *for every* $T > 0$ *with* $a(0) = \varphi(0)$, $\dot{a}(0) + v(0,0)\varphi'(0) = 0$ *and* $v \in C^1(\Re_+ \times [0,1])$ *is a non-negative function (i.e.,* $v(t,x) \geq 0$ *for all* $t \geq 0$, $x \in [0,1]$*) which has Lipschitz derivatives on* $[0,T]\times[0,1]$ *for every* $T > 0$. *Assume that* $v(t,0) > 0$ *for all* $t \geq 0$. *Then the initial-boundary value problem (5.3), (5.4), (5.5) has a unique solution* $w \in C^1(\Re_+ \times [0,1])$, *which has Lipschitz derivatives on* $[0,T]\times[0,1]$ *for every* $T > 0$ *and satisfies the inequality*

$$\|w[t]\|_\infty \leq \max\left(\max_{0 \leq s \leq t}(|a(s)|), \|\varphi\|_\infty\right), \text{ for all } t \geq 0. \tag{5.6}$$

*Moreover, if there exists a constant* $v_{\min} > 0$ *such that* $v(t,x) \geq v_{\min}$ *for all* $t \geq 0$, $x \in [0,1]$ *and if* $a \equiv 0$ *then* $w(t,x) = 0$ *for all* $x \in [0,1]$ *and* $t \geq v_{\min}^{-1}$.

*5.II. Proof of Theorem 2.1*

**Proof of Theorem 2.1:** Let arbitrary $T > 0$ be given. We will apply the method of finite differences (used in the book [14]) in order to construct a solution on $[0,T]$ for the initial-boundary value problem (2.11), (2.12), (2.13).

Let $N > c^{-1}\mu$ be an integer and consider the parameterized (with parameter $N$) discrete-time system

$$\begin{aligned}
w_i((k+1)\delta) &= (1 - \lambda v_i(k\delta))w_i(k\delta) + \lambda v_i(k\delta)w_{i-1}(k\delta) & i = 1,\ldots,N \\
v_i((k+1)\delta) &= (1 - \lambda c)v_i(k\delta) + \lambda c\, v_{i+1}(k\delta) & i = 0,\ldots,N-1 \\
v_N((k+1)\delta) &= (1 - \mu\delta)v_N(k\delta) + \mu\delta\, g(v_N(k\delta), w_N(k\delta))
\end{aligned} \tag{5.7}$$

where $k = 0,\ldots,m-1$ is an integer (time of the discrete-time system), $\lambda := \dfrac{T}{1 + \left[T\max\left(\max_{0 \leq x \leq 1}(\varphi(x)), v_{\max}, c\right)\right]}$, $m := N\left(1 + \left[T\max\left(\max_{0 \leq x \leq 1}(\varphi(x)), v_{\max}, c\right)\right]\right)$,

$$w_0(k\delta) = a(k\delta, v_0(k\delta)), \text{ for } k = 1,\ldots,m \tag{5.8}$$

$$h := 1/N\ ,\ \delta := \lambda h \tag{5.9}$$

and initial condition

$$w_i(0) - \theta(ih) = v_i(0) - \varphi(ih) = 0 \qquad i = 0,\ldots,N \tag{5.10}$$

Notice that definition $\lambda := \dfrac{T}{1+\left[T\max\left(\max_{0\leq x\leq 1}(\varphi(x)), v_{\max}, c\right)\right]}$ guarantees that $\lambda > 0$ is sufficiently small so that

$$\lambda \max\left(\max_{0\leq x\leq 1}(\varphi(x)), v_{\max}, c\right) \leq 1. \tag{5.11}$$

Moreover, definitions (5.9) and $\lambda := \dfrac{T}{1+\left[T\max\left(\max_{0\leq x\leq 1}(\varphi(x)), v_{\max}, c\right)\right]}$,

$m := N\left(1+\left[T\max\left(\max_{0\leq x\leq 1}(\varphi(x)), v_{\max}, c\right)\right]\right)$ guarantee that $T = m\delta$.

We next prove by induction on $k$ that

$$0 \leq v_i(k\delta) \leq \max\left(\max_{0\leq x\leq 1}(\varphi(x)), v_{\max}\right), \text{ for all } i = 0,\ldots,N \text{ and } k = 0,\ldots,m. \tag{5.12}$$

Indeed, by virtue of (5.10) it follows that (5.12) holds for $k = 0$. Using (5.7) and (5.11) we are in a position to guarantee that $0 \leq v_i((k+1)\delta) \leq \max\left(\max_{0\leq x\leq 1}(\varphi(x)), v_{\max}\right)$, for all $i = 0,\ldots,N-1$, provided that (5.12) holds for certain $k = 0,\ldots,m-1$. Moreover, using (2.10), (5.7) and the fact that $N > c^{-1}\mu$ (which together with (5.11) and (5.9) implies that $\mu\delta \leq 1$) we can guarantee that $0 \leq v_N((k+1)\delta) \leq \max\left(\max_{0\leq x\leq 1}(\varphi(x)), v_{\max}\right)$, provided that (5.12) holds for certain $k = 0,\ldots,m-1$. Consequently, we get that $0 \leq v_i((k+1)\delta) \leq \max\left(\max_{0\leq x\leq 1}(\varphi(x)), v_{\max}\right)$, for all $i = 0,\ldots,N$, provided that (5.12) holds for certain $k = 0,\ldots,m-1$. Therefore, (5.12) holds.

Define

$$B_T := \max\left\{|a(t,v)| : t \in [0,T], 0 \leq v \leq \max\left(\max_{0\leq x\leq 1}(\varphi(x)), v_{\max}\right)\right\}. \tag{5.13}$$

We next prove that

$$|w_i(k\delta)| \leq \max\left(\|\theta\|_\infty, B_T\right), \text{ for all } i = 0,\ldots, N \text{ and } k = 0,\ldots,m. \tag{5.14}$$

Indeed, by virtue of (5.10) it follows that (5.14) holds for $k = 0$. Suppose that (5.14) holds for all $i = 0,\ldots,N$ and for certain $k = 0,\ldots,m-1$. Using (5.7), (5.8), (5.11), (5.12) and the triangle inequality we are in a position to guarantee that $|w_i((k+1)\delta)| \leq \max\left(\|\theta\|_\infty, B_T\right)$, for all $i = 1,\ldots,N$. Using (5.7), (5.12), definitions (5.13), (13) and the triangle inequality, we can also guarantee that $|w_0((k+1)\delta)| \leq \max\left(\|\theta\|_\infty, B_T\right)$. Thus, (5.14) holds.

We next prove that

$$v_i(k\delta) \geq \min\left(\min_{0\leq x\leq 1}(\varphi(x)), \min\{g(0,w) : |w| \leq \max(B_T, \|\theta\|_\infty)\}\right),$$
$$\text{for all } i = 0,\ldots,N \text{ and } k = 0,\ldots,m. \tag{5.15}$$

Indeed, by virtue of (5.10) it follows that (5.15) holds for $k = 0$. Suppose that (5.15) holds for all $i = 0,\ldots,N$ and for certain $k = 0,\ldots,m-1$. Using (5.7) and (5.11) we are in a position to guarantee that $v_i((k+1)\delta) \geq \min\left(\min_{0\leq x\leq 1}(\varphi(x)), \min\{g(0,w) : |w| \leq \max(B_T, \|\theta\|_\infty)\}\right)$, for all $i = 0,\ldots,N-1$. Moreover, using (2.10), (5.7), (5.14) and the fact that $N > c^{-1}\mu$ (which together with (5.11) and (5.9) implies that $\mu\delta \leq 1$) we can guarantee that $v_N((k+1)\delta) \geq \min\left(\min_{0\leq x\leq 1}(\varphi(x)), \min\{g(0,w) : |w| \leq \max(B_T, \|\theta\|_\infty)\}\right)$. Thus, (5.15) holds.

We define for $(t,x) \in [0,T]\times[0,1]$ and for every integer $N > c^{-1}\mu$ (recall that $h = N^{-1}$, $\delta = \lambda h$, $m\delta = T$):

$$w(k\delta, x; N) = (i+1-xN)w_i(k\delta) + (xN-i)w_{i+1}(k\delta)$$
$$v(k\delta, x; N) = (i+1-xN)v_i(k\delta) + (xN-i)v_{i+1}(k\delta)$$
with $i = [xN]$, for $x \in [0,1)$, $k = 0,...,m$, (5.16)

$w(k\delta, 1; N) = w_N(k\delta)$ and $v(k\delta, 1; N) = v_N(k\delta)$, for $k = 0,...,m$, (5.17)

$$w(t, x; N) = (k+1-\lambda^{-1}tN)w(k\delta, x; N) + (\lambda^{-1}tN-k)w((k+1)\delta, x; N)$$
$$v(t, x; N) = (k+1-\lambda^{-1}tN)v(k\delta, x; N) + (\lambda^{-1}tN-k)v((k+1)\delta, x; N)$$
with $k = [\lambda^{-1}tN]$ for $x \in [0,1]$, $t \in [0,T]$. (5.18)

It follows from (5.12), (5.14), (5.15) and definitions (5.16), (5.17), (5.18) that the following inequalities hold for all $(t, x) \in [0,T] \times [0,1]$ and for every integer $N > c^{-1}\mu$:

$$\min\left(\min_{0 \leq x \leq 1}(\varphi(x)), \min\left\{g(0, w) : |w| \leq \max(B_T, \|\theta\|_\infty)\right\}\right) \leq v(t, x; N) \leq \max\left(\max_{0 \leq x \leq 1}(\varphi(x)), v_{\max}\right) \quad (5.19)$$

$$|w(t, x; N)| \leq \max(\|\theta\|_\infty, B_T) \quad (5.20)$$

Since the rest of proof is long, we need to describe the major steps in the proof.

Step 1: We first show that there exists a constant $L := L(T, \theta, \varphi, a) > 0$ such that for every integer $N > c^{-1}\mu$ the functions $w(\cdot; N)$, $v(\cdot; N)$ are Lipschitz on $[0,T] \times [0,1]$ with Lipschitz constant $L$. This step is very important because it allows the application of Arzela-Ascoli theorem. More specifically, it follows from (5.19), (5.20) that the sequences of functions $\{w(\cdot; N)\}_{N=N^*}^\infty$, $\{v(\cdot; N)\}_{N=N^*}^\infty$ with $N^* = [c^{-1}\mu] + 1$, are uniformly bounded and equicontinuous. Therefore, compactness of $[0,T] \times [0,1]$ and the Arzela-Ascoli theorem implies that there exist Lipschitz functions $w : [0,T] \times [0,1] \to \Re$, $v : [0,T] \times [0,1] \to \Re$ and subsequences $\{w(\cdot; q)\}_{q=N^*}^\infty \subset \{w(\cdot; N)\}_{N=N^*}^\infty$, $\{v(\cdot; q)\}_{q=N^*}^\infty \subset \{v(\cdot; N)\}_{N=N^*}^\infty$, which converge uniformly on $[0,T] \times [0,1]$ to $w$ and $v$, respectively. Moreover, the functions $w$ and $v$ are Lipschitz on $[0,T] \times [0,1]$ with same Lipschitz constant $L$ and satisfy the same bounds with $w(\cdot; N)$ and $v(\cdot; N)$, i.e., for all $(t, x) \in [0,T] \times [0,1]$ it holds that

$$\min\left(\min_{0 \leq x \leq 1}(\varphi(x)), \min\left\{g(0, w) : |w| \leq \max(B_T, \|\theta\|_\infty)\right\}\right) \leq v(t, x) \leq \max\left(\max_{0 \leq x \leq 1}(\varphi(x)), v_{\max}\right) \quad (5.21)$$

$$|w(t, x)| \leq \max(\|\theta\|_\infty, B_T) \quad (5.22)$$

Using the fact that $w$ and $v$ are Lipschitz on $[0,T] \times [0,1]$, (5.8), (5.10) and the fact that $\{w(\cdot; q)\}_{q=N^*}^\infty$, $\{v(\cdot; q)\}_{q=N^*}^\infty$ converge uniformly on $[0,T] \times [0,1]$ to $w$ and $v$, we can conclude that (2.13) holds and $w(t, 0) = a(t, v(t, 0))$ for all $t \in [0,T]$.

We remark that in what follows the convergent subsequences $\{w(\cdot; q)\}_{q=N^*}^\infty \subset \{w(\cdot; N)\}_{N=N^*}^\infty$, $\{v(\cdot; q)\}_{q=N^*}^\infty \subset \{v(\cdot; N)\}_{N=N^*}^\infty$, will be denoted by $\{w(\cdot; N)\}_{N=N^*}^\infty$, $\{v(\cdot; N)\}_{N=N^*}^\infty$.

Step 2: We define the function $\xi : [0,T] \to \Re$ by means of the equations

$$\dot{\xi}(t) = -\mu\left(\xi(t) - g(\xi(t), w(t,1))\right) = 0, \text{ for } t \in [0,T] \quad (5.23)$$

$$\xi(0) = \varphi(1) \quad (5.24)$$

and we show that $\xi(t) = v(t, 1)$ for all $t \in [0,T]$. Notice that $\xi \in W^{2,\infty}([0,T])$.

Step 3: We define the function $\tilde{v}:[0,T]\times[0,1]\to \Re$ by means of the formula

$$\tilde{v}(t,x) = \begin{cases} \varphi(x+ct) & if \quad x+ct \leq 1 \\ \xi\left(t - c^{-1}(1-x)\right) & if \quad x+ct > 1 \end{cases} \quad (5.25)$$

Due to the facts that $\xi \in W^{2,\infty}([0,T])$, $\varphi \in W^{2,\infty}([0,1])$ and since the compatibility conditions (5.24), $\varphi'(1) = -c^{-1}\mu(\varphi(1) - g(\varphi(1),\theta(1)))$ hold, it follows that $\tilde{v} \in C^1([0,T]\times[0,1])$ has Lipschitz derivatives satisfying $\frac{\partial \tilde{v}}{\partial t}(t,x) - c\frac{\partial \tilde{v}}{\partial x}(t,x) = 0$ for $(t,x) \in [0,T]\times[0,1]$, $\tilde{v}(0,x) = \varphi(x)$ for $x \in [0,1]$ and $\tilde{v}(t,1) = \xi(t) = v(t,1)$ for $t \in [0,T]$. We show that $\tilde{v}(t,x) = v(t,x)$ for all $(t,x) \in [0,T]\times[0,1]$. Thus, it follows from (5.23), (5.25) that the function $v$ is of class $C^1([0,T]\times[0,1])$ with Lipschitz derivatives and satisfies the equations $\frac{\partial v}{\partial t}(t,x) - c\frac{\partial v}{\partial x}(t,x) = 0$ for $(t,x) \in [0,T]\times[0,1]$, $v(0,x) = \varphi(x)$ for $x \in [0,1]$ and $\frac{\partial v}{\partial t}(t,1) = -\mu(v(t,1) - g(v(t,1),w(t,1)))$ for $t \in [0,T]$.

Step 4: Proposition 5.2 implies that there exists a unique $C^1$ solution $\tilde{w}:[0,T]\times[0,1]\to \Re$ of the initial-boundary value problem

$$\frac{\partial \tilde{w}}{\partial t}(t,x) + v(t,x)\frac{\partial \tilde{w}}{\partial x}(t,x) = 0, \text{ for all } (t,x) \in [0,T]\times[0,1] \quad (5.26)$$

$$\tilde{w}(t,0) = a(t,v(t,0)), \text{ for all } t \in [0,T] \quad (5.27)$$

$$\tilde{w}(0,x) = \theta(x), \text{ for all } x \in [0,1] \quad (5.28)$$

Moreover, $\tilde{w}$ has Lipschitz derivatives. We show that $\tilde{w}(t,x) = w(t,x)$ for all $(t,x) \in [0,T]\times[0,1]$. It thus follows that the functions $w,v$ are of class $C^1([0,T]\times[0,1])$ with Lipschitz derivatives and satisfy the equations (2.11), (2.12), (2.13) on $[0,T]\times[0,1]$. The fact that the solution satisfies estimates (2.14), (2.15) is a consequence of (5.21), (5.22) and the fact that $T>0$ is arbitrary.

Step 5: Finally, we prove that the constructed solution is unique.

We start with the proofs of all steps.

Step 1: Define for every $i = 0,...,N-1$ and $k = 0,...,m$:

$$\begin{aligned} y_i(k\delta) &= h^{-1}(w_{i+1}(k\delta) - w_i(k\delta)) \\ p_i(k\delta) &= h^{-1}(v_{i+1}(k\delta) - v_i(k\delta)) \end{aligned} \quad (5.29)$$

Using (5.7), (5.8), (5.9), (5.10), the fact that $f(0) = a(0,\varphi(0))$, we are in a position to verify that the following equations hold for all $k = 0,...,m-1$:

$$y_i((k+1)\delta) = (1-\lambda v_{i+1}(k\delta))y_i(k\delta) + \lambda v_{i+1}(k\delta) y_{i-1}(k\delta) - \delta p_i(k\delta) y_{i-1}(k\delta), \text{ for } i=1,...,N-1 \quad (5.30)$$

$$y_0((k+1)\delta) = (1-\lambda v_1(k\delta))y_0(k\delta) - \lambda \frac{a((k+1)\delta,v_0(k\delta)) - a(k\delta,v_0(k\delta))}{\delta} \\ - \frac{a((k+1)\delta,(1-\lambda c)v_0(k\delta) + \lambda c v_1(k\delta)) - a((k+1)\delta,v_0(k\delta))}{h} \quad (5.31)$$

$$p_i((k+1)\delta) = (1-\lambda c)p_i(k\delta) + \lambda c \, p_{i+1}(k\delta), \text{ for } i=0,1,...,N-2 \quad (5.32)$$

$$p_{N-1}((k+1)\delta) = (1-\lambda c)p_{N-1}(k\delta) + \mu\lambda(g(v_N(k\delta),w_N(k\delta)) - v_N(k\delta)) \quad (5.33)$$

Using (2.10), (5.11), (5.12), (5.15) and (5.32), (5.33), we get for all $k = 0,...,m-1$:

$$\max_{i=0,...N-1}\left(|p_i((k+1)\delta)|\right) \leq \max\left(\max_{i=0,...N-1}\left(|p_i(k\delta)|\right), (1-\lambda c)\max_{i=0,...N-1}\left(|p_i(k\delta)|\right) + \mu\lambda\overline{v}_{\max}\right) \quad (5.34)$$

where $\overline{v}_{\max} := \max\left(\max_{0\leq x\leq 1}(\varphi(x)), v_{\max}\right)$. Therefore, we obtain from (5.34) and Lemma 5.1 the estimate:

$$\max_{i=0,...N-1}\left(|p_i(k\delta)|\right) \leq \max_{i=0,...N-1}\left(|p_i(0)|\right) + c^{-1}\mu\overline{v}_{\max}, \text{ for all } k = 0,...,m \quad (5.35)$$

Definition (5.29) in conjunction with (5.10) implies that $|p_i(0)| \leq \|\varphi'\|_\infty$ for all $i = 0,...,N-1$. Consequently, we get from (5.35) that

$$\max_{i=0,...N-1}\left(|p_i(k\delta)|\right) \leq P := \|\varphi'\|_\infty + c^{-1}\mu\overline{v}_{\max}, \text{ for all } k = 0,...,m \quad (5.36)$$

Using (5.11), (5.12), (5.15), (5.30), (5.31), (5.36), we get for all $k = 0,...,m-1$:

$$\max_{i=0,...N-1}\left(|y_i((k+1)\delta)|\right) \leq \max\left((1+\delta P)\max_{i=0,...N-1}\left(|y_i(k\delta)|\right), (1-\lambda v_{\min})\max_{i=0,...N-1}\left(|y_i(k\delta)|\right) + \lambda R(1+cP)\right) \quad (5.37)$$

where $v_{\min} := \min\left(\min_{0\leq x\leq 1}(\varphi(x)), \min\{g(0,w):|w|\leq \max(B_T,\|\theta\|_\infty)\}\right)$, $\overline{v}_{\max} := \max\left(\max_{0\leq x\leq 1}(\varphi(x)), v_{\max}\right)$ and
$R := \max\left\{\left|\frac{\partial a}{\partial t}(t,v)\right| + \left|\frac{\partial a}{\partial v}(t,v)\right|: 0\leq v\leq \overline{v}_{\max}, t\in[0,T]\right\}$. Using (5.37), Lemma 5.1 and the fact that $m\delta = T$, we get for all $k = 0,...,m$:

$$\max_{i=0,...N-1}\left(|y_i(k\delta)|\right) \leq \exp(PT)\left(\max_{i=0,...N-1}\left(|y_i(0)|\right) + R\frac{1+cP}{v_{\min}}\right) \quad (5.38)$$

Definition (5.29) in conjunction with (5.10) implies that $|y_i(0)| \leq \|\theta'\|_\infty$ for all $i = 0,...,N-1$. Consequently, we get from (5.35) that

$$\max_{i=0,...N-1}\left(|y_i(k\delta)|\right) \leq Y := \exp(PT)\left(\|\theta'\|_\infty + R\frac{1+cP}{v_{\min}}\right), \text{ for all } k = 0,...,m \quad (5.39)$$

It follows from (5.29), (5.36) and (5.39) that the following inequalities hold for all $i,j \in \{0,...,N\}$ and $k = 0,...,m$:

$$|w_i(k\delta) - w_j(k\delta)| \leq h|i-j|Y \text{ and } |v_i(k\delta) - v_j(k\delta)| \leq h|i-j|P \quad (5.40)$$

Notice that $P, Y$, defined in (5.36) and (5.39), respectively, depend only on $T, \theta, \varphi$ and $a$.

Next define for every $i = 0,...,N$ and $k = 0,...,m-1$:

$$\begin{aligned}\zeta_i(k\delta) &= \delta^{-1}(w_i((k+1)\delta) - w_i(k\delta)) \\ \eta_i(k\delta) &= \delta^{-1}(v_i((k+1)\delta) - v_i(k\delta))\end{aligned} \quad (5.41)$$

It follows from (5.7) and definitions (5.9), (5.29), (5.41) that the following equalities hold for $k = 0,...,m-1$:

$$\begin{aligned}\zeta_i(k\delta) &= -v_i(k\delta)y_{i-1}(k\delta) & i &= 1,...,N \\ \eta_i(k\delta) &= c\,p_i(k\delta) & i &= 0,...,N-1 \\ \eta_N(k\delta) &= \mu(g(v_N(k\delta), w_N(k\delta)) - v_N(k\delta))\end{aligned} \quad (5.42)$$

Using (5.42), (2.10), (5.12) and (5.36) we get:

$$\max_{i=0,...N}\left(|\eta_i(k\delta)|\right) \leq cP + \mu\overline{v}_{\max}, \text{ for all } k = 0,...,m-1 \quad (5.43)$$

where $\overline{v}_{\max} := \max\left(\max_{0\leq x\leq 1}(\varphi(x)), v_{\max}\right)$. It follows from (5.41), (5.7), (5.8), (5.10) and the fact that $f(0) = a(0, \varphi(0))$ that the following equalities hold for $k = 0,...,m-1$:

$$\zeta_0(k\delta) = \frac{a((k+1)\delta, v_0(k\delta)) - a(k\delta, v_0(k\delta))}{\delta} + \frac{a((k+1)\delta, (1-\lambda c)v_0(k\delta) + \lambda cv_1(k\delta)) - a((k+1)\delta, v_0(k\delta))}{\delta} \quad (5.44)$$

Equalities (5.42), (5.44) in conjunction with (5.12), (5.36), (5.39) and definitions (5.9), (5.29) imply that

$$\max_{i=0,\dots,N}\left(|\zeta_i(k\delta)|\right) \leq R(1+cP) + \bar{v}_{\max} Y, \text{ for all } k = 0,\dots,m-1 \quad (5.45)$$

where $\bar{v}_{\max} := \max\left(\max_{0\leq x\leq 1}(\varphi(x)), v_{\max}\right)$ and $R := \max\left\{\left|\frac{\partial a}{\partial t}(t,v)\right| + \left|\frac{\partial a}{\partial v}(t,v)\right| : 0 \leq v \leq \bar{v}_{\max}, t \in [0,T]\right\}$.

It follows from (5.41), (5.43) and (5.45) that the following inequalities hold for all $i = 0,\dots,N$ and $k,l \in \{0,\dots,m\}$:

$$|w_i(k\delta) - w_i(l\delta)| \leq \delta|k-l|\left(R(1+cP) + \bar{v}_{\max} Y\right)$$
$$|v_i(k\delta) - v_i(l\delta)| \leq \delta|k-l|\left(cP + \mu \bar{v}_{\max}\right) \quad (5.46)$$

It follows from definitions (5.16), (5.17), (5.18) and inequalities (5.40), (5.46) that there exists a constant $L := L(T,\theta,\varphi,a) > 0$ such that for every integer $N > c^{-1}\mu$ the following Lipschitz inequalities hold:

$$|w(t,x;N) - w(\tau,z;N)| + |v(t,x;N) - v(\tau,z;N)| \leq L(|x-z| + |t-\tau|), \text{ for all } x,z \in [0,1], t,\tau \in [0,T] \quad (5.47)$$

This completes Step 1 of the proof.

<u>Step 2:</u> We define $\xi:[0,T] \to \Re$ by means of (5.23), (5.24). The fact that $\xi$ can be defined on the whole $[0,T]$ is a consequence of the formula $\xi(t) = \exp(-\mu t)\varphi(1) + \mu \int_0^t \exp(-\mu(t-s))g(\xi(s),w(s,1))ds$, which together with (2.10) and (5.22) implies the inequality $\min\left(\varphi(1), \min\{g(0,w) : |w| \leq \max(B_T, \|\theta\|_\infty)\}\right) \leq \xi(t) \leq \max(\varphi(1), v_{\max})$ for all $t \in [0,T]$. Pick any integer $N > c^{-1}\mu$. It follows from (2.14) and (5.23) that the following inequality holds for $k = 0,\dots,m-1$:

$$\left|\xi((k+1)\delta) - (1-\mu\delta)\xi(k\delta) - \mu\delta g(\xi(k\delta), w(k\delta,1))\right| \leq \mu \frac{\delta^2}{2}\left(\|\dot{\xi}\|(1+G) + GL\right) \quad (5.48)$$

where $m := N\left(1 + \left[T\max\left(\max_{0\leq x\leq 1}(\varphi(x)), v_{\max}, c\right)\right]\right)$, $\delta = \lambda/N$, $\lambda := \frac{T}{1 + \left[T\max\left(\max_{0\leq x\leq 1}(\varphi(x)), v_{\max}, c\right)\right]}$,

$G := \max\left\{\left|\frac{\partial g}{\partial v}(v,w)\right| + \left|\frac{\partial g}{\partial w}(v,w)\right| : (v,w) \in S\right\}$, $L$ is the Lipschitz constant of $w$, $S := \left\{(v,w) \in \Re_+ \times \Re : v \leq \max(\varphi(1), v_{\max}), |w| \leq \max(B_T, \|f\|_\infty)\right\}$ and $\|\dot{\xi}\| = \max_{0\leq t\leq T}\left(|\dot{\xi}(t)|\right)$. Using (5.7), (5.48) and defining $e_N(k\delta) := \xi(k\delta) - v_N(k\delta)$ we get for $k = 0,\dots,m-1$:

$$|e_N((k+1)\delta)| \leq (1 + \mu\delta G)|e_N(k\delta)| + \mu\delta G|w(k\delta,1) - w(k\delta,1;N)| + \mu\frac{\delta^2}{2}\left(\|\dot{\xi}\|(1+G) + LG\right) \quad (5.49)$$

Inequality (5.49) in conjunction with Lemma 5.1 and the facts that $T = m\delta$, $e_N(0) = 0$ (a consequence of (5.24) and (2.13)) implies that the following estimate holds:

$$|e_N(k\delta)| \leq \mu cT \exp(\mu cTG)\left(\max_{0\leq t\leq T}(|w(t,1) - w(t,1;N)|) + \frac{\delta}{2}\left(\|\dot{\xi}\|(1+G) + LG\right)\right), \text{ for } k = 0,\dots,m \quad (5.50)$$

Pick any $t \in [0,T]$ and set $k = [t\delta^{-1}]$. Then we get

$$|\xi(t) - v(t,1)| \leq |\xi(t) - v(t,1) - \xi(k\delta) + v(k\delta,1)| + |\xi(k\delta) - v(k\delta,1)|$$
$$\leq (L + \|\dot{\xi}\|)(t - k\delta) + |\xi(k\delta) - v_N(k\delta)| + |v(k\delta,1) - v(k\delta,1;N)|$$
$$\leq (L + \|\dot{\xi}\|)\delta + |e_N(k\delta)| + \max_{0\leq s\leq T}|v(s,1) - v(s,1;N)|$$

where $L$ is the Lipschitz constant of $v$. In the above derivation, we have used the fact that $t - k\delta \leq \delta$. Since $\{w(\cdot;N)\}_{N=N^*}^\infty$, $\{v(\cdot;N)\}_{N=N^*}^\infty$ converge uniformly to $w$ and $v$ as $N \to +\infty$, the above inequality in conjunction with (5.50) shows that $\xi(t) = v(t,1)$ for all $t \in [0,T]$. This completes Step 2 of the proof.

**Step 3:** Pick any integer $N > c^{-1}\mu$. Define $m := N\left(1 + \left[T\max\left(\max_{0 \leq x \leq 1}(\varphi(x)), v_{\max}, c\right)\right]\right)$, $h = 1/N$, $\delta = \lambda/N$,

$\lambda := \dfrac{T}{1 + \left[T\max\left(\max_{0 \leq x \leq 1}(\varphi(x)), v_{\max}, c\right)\right]}$ and notice that since $\tilde{v} \in C^1([0,T] \times [0,1])$ satisfies

$\dfrac{\partial \tilde{v}}{\partial t}(t,x) - c\dfrac{\partial \tilde{v}}{\partial x}(t,x) = 0$ for $(t,x) \in [0,T] \times [0,1]$, $\tilde{v}(0,x) = \varphi(x)$ for $x \in [0,1]$ and $\tilde{v}(t,1) = \xi(t) = v(t,1)$ for $t \in [0,T]$ then we get for all $k = 0,\ldots,m-1$ and $i = 0,\ldots,N-1$:

$$\tilde{v}((k+1)\delta, ih) = (1 - \lambda c)\tilde{v}(k\delta, ih) + c\lambda \tilde{v}(k\delta, (i+1)h) + c\, err(k,i) \tag{5.51}$$

where

$$err(k,i) = \int_{k\delta}^{(k+1)\delta}\left(\dfrac{\partial \tilde{v}}{\partial x}(t, ih) - \dfrac{\partial \tilde{v}}{\partial x}(k\delta, ih)\right)dt - \lambda \int_{ih}^{(i+1)h}\left(\dfrac{\partial \tilde{v}}{\partial x}(k\delta, x) - \dfrac{\partial \tilde{v}}{\partial x}(k\delta, ih)\right)dx \tag{5.52}$$

Defining $e_i^v(k\delta) := \tilde{v}(k\delta, ih) - v(k\delta, ih; N)$ for $k = 0,\ldots,m$ and $i = 0,\ldots,N$, we get from (5.7), (5.10), (5.16), (5.17), (5.51) and the facts that $\tilde{v}(0,x) = \varphi(x)$ for $x \in [0,1]$ and $\tilde{v}(t,1) = \xi(t) = v(t,1)$ for $t \in [0,T]$:

$$e_i^v((k+1)\delta) = (1 - \lambda c)e_i^v(k\delta, ih) + c\lambda e_{i+1}^v(k\delta, (i+1)h) + c\, err(k,i), \text{ for } k = 0,\ldots,m-1, \; i = 0,\ldots,N-1 \tag{5.53}$$

$$e_i^v(0) = 0, \text{ for } i = 0,\ldots,N \tag{5.54}$$

$$e_N^v(k\delta) = 0, \text{ for } k = 0,\ldots,m \tag{5.55}$$

Using (5.9), (5.11), (5.52), (5.53), (5.55), we get:

$$\max_{i=0,\ldots,N}\left(\left|e_i^v((k+1)\delta)\right|\right) \leq \max_{i=0,\ldots,N}\left(\left|e_i^v(k\delta)\right|\right) + 2c\delta\, G(N), \text{ for } k = 0,\ldots,m-1 \tag{5.56}$$

where

$$G(N) := \max\left\{\left|\dfrac{\partial \tilde{v}}{\partial x}(t,x) - \dfrac{\partial \tilde{v}}{\partial x}(\tau,z)\right| : (t,x) \in [0,T] \times [0,1], (\tau,z) \in [0,T] \times [0,1], |t-\tau| + |x-z| \leq (1+\lambda)N^{-1}\right\} \tag{5.57}$$

Using Lemma 5.1, in conjunction with (5.54), (5.56) and the fact that $T = m\delta$, we get $\max_{i=0,\ldots,N}\left(\left|e_i^v(k\delta)\right|\right) \leq 2cT\, G(N)$. Next, pick any $(t,x) \in [0,T] \times [0,1]$ and set $k = [t\delta^{-1}]$, $i = [xh^{-1}]$. Using the inequality $|\tilde{v}(t,x) - v(t,x)| \leq |\tilde{v}(t,x) - \tilde{v}(k\delta, ih)| + |e_i^v(k\delta)| + |v(k\delta, ih; N) - v(k\delta, ih)| + |v(t,x) - v(k\delta, ih)|$, we are in a position to conclude (using the fact that $t - k\delta \leq \delta$ and $x - ih \leq h$) that

$$|\tilde{v}(t,x) - v(t,x)| \leq 2L(1+\lambda)N^{-1} + 2cTG(N) + \max_{(\tau,z) \in [0,T] \times [0,1]}\left(|v(\tau,z;N) - v(\tau,z)|\right) \tag{5.58}$$

where $L$ is the Lipschitz constant of $v$ and $\tilde{v}$. Definition (5.57), the fact that $\tilde{v} \in C^1([0,T] \times [0,1])$ has Lipschitz derivatives on $[0,T] \times [0,1]$ implies that $\lim_{N \to +\infty}(G(N)) = 0$. Moreover, since $\{v(\cdot; N)\}_{N=N^*}^{\infty}$ converges uniformly to $v$ as $N \to +\infty$, we get from (5.58) that $\tilde{v}(t,x) = v(t,x)$ for all $(t,x) \in [0,T] \times [0,1]$. This completes Step 3 of the proof.

**Step 4:** Pick any integer $N > c^{-1}\mu$. Define $m := N\left(1 + \left[T\max\left(\max_{0 \leq x \leq 1}(\varphi(x)), v_{\max}, c\right)\right]\right)$, $h = 1/N$, $\delta = \lambda/N$,

$\lambda := \dfrac{T}{1 + \left[T\max\left(\max_{0 \leq x \leq 1}(\varphi(x)), v_{\max}, c\right)\right]}$ and notice that since $\tilde{w} \in C^1([0,T] \times [0,1])$ satisfies (5.26) for all $(t,x) \in [0,T] \times [0,1]$, we get for all $k = 0,\ldots,m-1$ and $i = 1,\ldots,N$:

$$\tilde{w}((k+1)\delta, ih) = (1 - \lambda v(k\delta, ih))\tilde{w}(k\delta, ih) + \lambda v(k\delta, ih)\tilde{w}(k\delta, (i-1)h) + Err(k,i) \tag{5.59}$$

where

$$Err(k,i) = \lambda v(k\delta, ih) \int_{(i-1)h}^{ih} \left( \frac{\partial \tilde{w}}{\partial x}(k\delta, x) - \frac{\partial \tilde{w}}{\partial x}(k\delta, ih) \right) dx - \int_{k\delta}^{(k+1)\delta} \left( v(t, ih) \frac{\partial \tilde{w}}{\partial x}(t, ih) - v(k\delta, ih) \frac{\partial \tilde{w}}{\partial x}(k\delta, ih) \right) dt \quad (5.60)$$

Defining $e_i^w(k\delta) := \tilde{w}(k\delta, ih) - w(k\delta, ih; N)$ for $k = 0,...,m$ and $i = 0,...,N$, we get from (5.7), (5.9), (5.10), (5.16), (5.17), (5.27), (5.28), (5.29), (5.59):

$$e_i^w((k+1)\delta) = (1 - \lambda v(k\delta, ih)) e_i^w(k\delta) - \delta (v(k\delta, ih) - v(k\delta, ih; N)) y_{i-1}(k\delta) + \lambda v(k\delta, ih) e_{i-1}^w(k\delta) + Err(k,i) ,$$
$$\text{for } k = 0,...,m-1 , \ i = 1,...,N \quad (5.61)$$

$$e_i^w(0) = 0 , \text{ for } i = 0,...,N \quad (5.62)$$

$$e_0^w(k\delta) = 0 , \text{ for } k = 0,...,m \quad (5.63)$$

Using (2.15), (5.9), (5.11), (5.39), (5.60), (5.61), (5.63), we get:

$$\max_{i=0,...N} \left( \left| e_i^w((k+1)\delta) \right| \right) \leq \max_{i=0,...N} \left( \left| e_i^w(k\delta) \right| \right) + \delta \tilde{G}(N) , \text{ for } k = 0,...,m-1 \quad (5.64)$$

where

$$\tilde{G}(N) := \bar{v}_{\max} \max \left\{ \left| \frac{\partial \tilde{w}}{\partial x}(t,x) - \frac{\partial \tilde{w}}{\partial x}(\tau, z) \right| : (t,x) \in [0,T] \times [0,1], (\tau, z) \in [0,T] \times [0,1], |t-\tau| + |x-z| \leq (1+\lambda)N^{-1} \right\}$$
$$+ \lambda N^{-1} L \max \left\{ \left| \frac{\partial \tilde{w}}{\partial x}(t,x) \right| : (t,x) \in [0,T] \times [0,1] \right\} + Y \max_{(t,x) \in [0,T] \times [0,1]} \left( |v(t,x) - v(t,x;N)| \right) \quad (5.65)$$

and $\bar{v}_{\max} = \max \left( \max_{0 \leq x \leq 1} (\varphi(x)), v_{\max} \right)$. Using Lemma 5.1, in conjunction with (5.62), (5.65) and the fact that $T = m\delta$, we get $\max_{i=0,...N} \left( \left| e_i^w(k\delta) \right| \right) \leq T \tilde{G}(N)$. Next, pick any $(t,x) \in [0,T] \times [0,1]$ and set $k = [t\delta^{-1}]$, $i = [xh^{-1}]$. Using the inequality $|\tilde{w}(t,x) - w(t,x)| \leq |\tilde{w}(t,x) - \tilde{w}(k\delta, ih)| + |e_i^w(k\delta)| + |w(k\delta, ih; N) - w(k\delta, ih)| + |w(t,x) - w(k\delta, ih)|$, we are in a position to conclude (using the fact that $t - k\delta \leq \delta$ and $x - ih \leq h$) that

$$|\tilde{v}(t,x) - v(t,x)| \leq 2L(1+\lambda)N^{-1} + T\tilde{G}(N) + \max_{(\tau, z) \in [0,T] \times [0,1]} \left( |w(\tau, z; N) - w(\tau, z)| \right) \quad (5.66)$$

where $L$ is the Lipschitz constant of $w$ and $\tilde{w}$. Definition (5.65), in conjunction with the fact that $\tilde{w} \in C^1([0,T] \times [0,1])$ has Lipschitz derivatives on $[0,T] \times [0,1]$ implies that $\lim_{N \to +\infty} (\tilde{G}(N)) = 0$. Moreover, since $\{w(\cdot; N)\}_{N=N^*}^{\infty}$ converges uniformly to $w$ as $N \to +\infty$, we get from (5.66) that $\tilde{w}(t,x) = w(t,x)$ for all $(t,x) \in [0,T] \times [0,1]$. This completes Step 4 of the proof.

Step 5: Uniqueness follows from a contradiction argument. Suppose that there exist two solutions $(w,v) \in C^1(\Re_+ \times [0,1]) \times C^1(\Re_+ \times [0,1])$ and $(\bar{w}, \bar{v}) \in C^1(\Re_+ \times [0,1]) \times C^1(\Re_+ \times [0,1])$ of the initial-boundary value problem (2.10), (2.11), (2.12). It then follows that the functions $e_w = \bar{w} - w$, $e_v = \bar{v} - v$ satisfy the following equations:

$$\frac{\partial e_w}{\partial t}(t,x) + v(t,x) \frac{\partial e_w}{\partial x}(t,x) + e_v(t,x) \frac{\partial \bar{w}}{\partial x}(t,x) = \frac{\partial e_v}{\partial t}(t,x) - c \frac{\partial e_v}{\partial x}(t,x) = 0 , \text{ for all } (t,x) \in \Re_+ \times [0,1] \quad (5.67)$$

$$e_w(t,0) = a(t, v(t,0) + e_v(t,0)) - a(t, v(t,0)) , \text{ for all } t \geq 0 \quad (5.68)$$

$$\frac{\partial e_v}{\partial t}(t,1) = -\mu \left( e_v(t,1) + g(v(t,1), w(t,1)) - g(\bar{v}(t,1) + e_v(t,1), \bar{w}(t,1) + e_w(t,1)) \right) , \text{ for all } t \geq 0 \quad (5.69)$$

$$e_w(0,x) = e_v(0,x) = 0 , \text{ for all } x \in [0,1] \quad (5.70)$$

Let $T > 0$ be given and let $S \subseteq \Re_+ \times \Re$ be a compact set that contains both solutions on $[0,T]$, i.e., $(v(t,x), w(t,x)) \in S$ and $(\bar{v}(t,x), \bar{w}(t,x)) \in S$ for all $(t,x) \in [0,T] \times [0,1]$. Let $M \geq 1$ be a constant that satisfies $Mc \geq Q^2 \bar{v}_{max}$, where $Q := \max\left\{\left|\frac{\partial g}{\partial v}(v,w)\right| + \left|\frac{\partial g}{\partial w}(v,w)\right| : (v,w) \in S\right\} + \max\left\{\left|\frac{\partial a}{\partial v}(t,v)\right| : t \in [0,T], (v,w) \in S\right\}$ and

$\bar{v}_{max} = \max\left(\max_{0 \leq x \leq 1}(\varphi(x)), v_{max}\right)$. Define the functional:

$$V(t) = \frac{1}{2}\int_0^1 e_w^2(t,x)dx + \frac{M}{2}\int_0^1 e_v^2(t,x)dx + \frac{1}{2}e_v^2(t,1) \tag{5.71}$$

We next estimate the time derivative of $V(t)$ on $[0,T]$. Using (5.67), (5.68), (5.69), the fact that $|g(v(t,1),w(t,1)) - g(\bar{v}(t,1) + e_v(t,1), \bar{w}(t,1) + e_w(t,1))| \leq Q(|e_v(t,1)| + |e_w(t,1)|)$ and the fact that $|e_v(t,x)e_w(t,x)| \leq \frac{1}{2}e_v^2(t,x) + \frac{1}{2}e_w^2(t,x)$, we get:

$$\dot{V}(t) = -\int_0^1 v(t,x)e_w(t,x)\frac{\partial e_w}{\partial x}(t,x)dx - \int_0^1 e_v(t,x)e_w(t,x)\frac{\partial \bar{w}}{\partial x}(t,x)dx + Mc\int_0^1 e_v(t,x)\frac{\partial e_v}{\partial x}(t,x)dx + e_v(t,1)\frac{\partial e_v}{\partial t}(t,1)$$

$$\leq -\frac{1}{2}v(t,1)e_w^2(t,1) + \frac{1}{2}v(t,0)e_w^2(t,0) + \frac{1}{2}\|\bar{w}_x\|\left(\int_0^1 e_w^2(t,x)dx + \int_0^1 e_v^2(t,x)dx\right) + \frac{1}{2}\int_0^1 e_w^2(t,x)\frac{\partial v}{\partial x}(t,x)dx$$

$$+ \left(\frac{Mc}{2} - \mu\right)e_v^2(t,1) - \frac{Mc}{2}e_v^2(t,0) + \mu Q|e_v(t,1)|(|e_v(t,1)| + |e_w(t,1)|)$$

where $\|\bar{w}_x\| := \max\left\{\left|\frac{\partial \bar{w}}{\partial x}(t,x)\right| : (t,x) \in [0,T] \times [0,1]\right\}$. The above inequality in conjunction with (2.15), definition $v_{min} := \min\left(\min_{0 \leq x \leq 1}(\varphi(x)), \min\{g(0,w) : |w| \leq \max(B_T, \|\theta\|_\infty)\}\right)$ and the fact that $\mu Q|e_v(t,1)||e_w(t,1)| \leq \frac{1}{2}v_{min}e_w^2(t,1) + \frac{Q^2\mu^2}{2v_{min}}e_v^2(t,1)$, gives:

$$\dot{V}(t) \leq \frac{1}{2}v(t,0)e_w^2(t,0) + \frac{1}{2}(\|\bar{w}_x\| + \|v_x\|)\left(\int_0^1 e_w^2(t,x)dx + \int_0^1 e_v^2(t,x)dx\right) + \left(\frac{Mc}{2} + Q + \frac{Q^2\mu^2}{2v_{min}}\right)e_v^2(t,1) - \frac{Mc}{2}e_v^2(t,0)$$

where $\|v_x\| := \max\left\{\left|\frac{\partial v}{\partial x}(t,x)\right| : (t,x) \in [0,T] \times [0,1]\right\}$. Since $|a(t, v(t,0) + e_v(t,0)) - a(t, v(t,0))| \leq Q|e_v(t,0)|$, we get from the above inequality, the fact that $v(t,0) \leq \bar{v}_{max}$, where $\bar{v}_{max} = \max\left(\max_{0 \leq x \leq 1}(\varphi(x)), v_{max}\right)$, the facts that $M \geq 1$ and $Mc \geq Q^2\bar{v}_{max}$:

$$\dot{V}(t) \leq \frac{1}{2}(\|\bar{w}_x\| + \|v_x\|)\left(\int_0^1 e_w^2(t,x)dx + M\int_0^1 e_v^2(t,x)dx\right) + \left(\frac{Mc}{2} + Q + \frac{Q^2\mu^2}{2v_{min}}\right)e_v^2(t,1)$$

The above inequality in conjunction with definition (5.71) shows that there exists a constant $K > 0$ such that $\dot{V}(t) \leq KV(t)$ for all $t \in [0,T]$. Gronwall's lemma in conjunction with (5.70), (5.71) (which give $V(0) = 0$) implies that $V(t) \equiv 0$ on $[0,T]$. Therefore, since $T > 0$ is arbitrary, we conclude that $\bar{w}(t,x) - w(t,x) = \bar{v}(t,x) - v(t,x) = 0$ for all $x \in [0,1]$ and $t \geq 0$. The proof is complete. ◁

*5.III. Proof of Theorem 3.1*

**Proof of Theorem 3.1:** Let arbitrary $(\rho_0, v_0) \in X$ be given, for which the equalities $\rho_0(0) = \rho_{eq}\frac{c + f(\rho_{eq})}{c + v_0(0)}$, $\rho_0'(0) = -\frac{-\rho_0(0)}{c + v_0(0)}v_0'(0)$ hold. The solution of the initial-boundary value problem (2.1), (2.2), (2.3), (2.6) with (3.2), (3.3) is constructed by applying the transformation (2.7), using (2.4) and the condition $\rho_{eq} \leq \frac{c}{c+f(\rho_{eq})}(\rho_{max}-\varepsilon)$. More specifically, we get the initial-boundary value problem

$$\frac{\partial w}{\partial t}(t,x)+v(t,x)\frac{\partial w}{\partial x}(t,x)=\frac{\partial v}{\partial t}(t,x)-c\frac{\partial v}{\partial x}(t,x)=0, \text{ for all } (t,x)\in\Re_+\times[0,1] \quad (5.72)$$

$$w(t,0)=\frac{\partial v}{\partial t}(t,1)+\mu\left(v(t,1)-f\left(\rho_{eq}\exp(w(t,1))\frac{c+f(\rho_{eq})}{c+v(t,1)}\right)\right)=0, \text{ for all } t\geq 0 \quad (5.73)$$

$$w(0,x)-w_0(x)=v(0,x)-v_0(x)=0, \text{ for all } x\in[0,1] \quad (5.74)$$

where

$$w_0(x)=\ln\left(\frac{\rho_0(x)(c+v_0(x))}{(c+f(\rho_{eq}))\rho_{eq}}\right), \text{ for all } x\in[0,1] \quad (5.75)$$

and $\rho(t,x)$ may be found by means of (2.7). Exploiting the fact that $f$ is non-increasing and Theorem 2.1, we conclude that the initial-boundary value problem (5.72), (5.73), (5.74) admits a unique solution $w,v\in C^1(\Re_+\times[0,1])$, which has Lipschitz derivatives on every compact $S\subset\Re_+\times[0,1]$ and satisfies the following inequalities for all $(t,x)\in\Re_+\times[0,1]$:

$$\|w[t]\|_\infty \leq \|w_0\|_\infty \quad (5.76)$$

$$\min\left(\min_{0\leq x\leq 1}(v_0(x)), f\left(\rho_{eq}\exp(\|w_0\|_\infty)\frac{c+f(\rho_{eq})}{c}\right)\right) \leq v(t,x) \leq \max\left(\max_{0\leq x\leq 1}(v_0(x)), f(0)\right) \quad (5.77)$$

Define for all $(t,x)\in\Re_+\times[0,1]$:

$$b(t,x):=\ln\left(\frac{v(t,x)}{f(\rho_{eq})}\right), \quad b_0(x):=\ln\left(\frac{v_0(x)}{f(\rho_{eq})}\right) \quad (5.78)$$

Notice that (5.77) in conjunction with definition (5.78) gives:

$$v_{\min}:=\min\left(f(\rho_{eq})\exp(-\|b_0\|_\infty), f\left(\rho_{eq}\exp(\|w_0\|_\infty)\frac{c+f(\rho_{eq})}{c}\right)\right) \leq v(t,x) \leq \max\left(f(\rho_{eq})\exp(\|b_0\|_\infty), f(0)\right) \quad (5.79)$$

Moreover, Proposition 5.2 implies that

$$w(t,x)=0 \text{ for all } x\in[0,1] \text{ and } t\geq v_{\min}^{-1} \quad (5.80)$$

Equations (5.72) and (5.73) imply that the following equations holds for all $(t,x)\in\Re_+\times[0,1]$:

$$v(t,x)=\begin{cases}v_0(x+ct) & \text{if } x+ct\leq 1 \\ \xi(t-c^{-1}(1-x)) & \text{if } x+ct>1\end{cases} \quad (5.81)$$

where $\xi:\Re_+\to\Re$ is the solution of the initial-value problem

$$\dot\xi(t)=-\mu\left(\xi(t)-f\left(\rho_{eq}\exp(w(t,1))\frac{c+f(\rho_{eq})}{c+\xi(t)}\right)\right) \quad (5.82)$$

$$\xi(0)=v_0(1). \quad (5.83)$$

Formula (5.81) implies the following estimate for every $\sigma>0$, $t\geq 0$:

$$\|b[t]\|_\infty \leq \exp(-\sigma(ct-1))\|b_0\|_\infty + \max_{\max(0,t-c^{-1})\leq s\leq t}\left(\left|\ln\left(\frac{\xi(s)}{f(\rho_{eq})}\right)\right|\right) \quad (5.84)$$

Using the transformation

$$\zeta(t) = \ln\left(\frac{\xi(t)}{f(\rho_{eq})}\right), \quad (5.85)$$

we get from (5.82):

$$\dot\zeta(t) = -\mu\left(1 - \frac{\exp(-\zeta(t))}{f(\rho_{eq})}f\left(\rho_{eq}\exp(w(t,1))\frac{c+f(\rho_{eq})}{c+f(\rho_{eq})\exp(\zeta(t))}\right)\right) \quad (5.86)$$

Inequality (3.1) implies that $0 \in \Re$ is a globally asymptotically stable equilibrium point for system (5.82) with $w(t,1) \equiv 0$. Consequently, it follows from (5.80) and Theorem 2.2 in [15] that there exists a function $P \in KL$ such that the following estimate holds:

$$\left|\ln\left(\frac{\xi(t)}{f(\rho_{eq})}\right)\right| \leq P\left(\left|\ln\left(\frac{\xi(v_{\min}^{-1})}{f(\rho_{eq})}\right)\right|, t - v_{\min}^{-1}\right), \text{ for all } t \geq v_{\min}^{-1} \quad (5.87)$$

Since $(\zeta(t), w(t,1)) \in \Re^2$ takes values in a compact set $S(w_0, v_0) \subset \Re^2$ for all $t \geq 0$ (recall definition (5.85), (5.76), (5.79) and (5.81)) and since $F(\zeta, w) := -\mu\left(1 - \frac{\exp(-\zeta)}{f(\rho_{eq})}f\left(\rho_{eq}\exp(w)\frac{c+f(\rho_{eq})}{c+f(\rho_{eq})\exp(\zeta)}\right)\right)$ is a $C^1$ mapping, there exists a non-decreasing function $L: \Re_+ \to \Re_+$ such that $|\dot\zeta(t)| \leq (|\zeta(t)| + |w(t,1)|)L(\|w_0\|_\infty + \|b_0\|_\infty)$. Using Gronwall's Lemma in conjunction with (5.76) and the previous inequality we get $|\zeta(t)| \leq \exp(Lv_{\min}^{-1})(|\zeta(0)| + Lv_{\min}^{-1}\|w_0\|_\infty)$, for all $t \in [0, v_{\min}^{-1}]$, where $L := L(\|w_0\|_\infty + \|b_0\|_\infty)$. Combining (5.87) with the previous estimate, definitions (5.78), (5.85) and (5.83), we get for all $t \geq 0$:

$$\left|\ln\left(\frac{\xi(t)}{f(\rho_{eq})}\right)\right| \leq P\left(\exp(Lv_{\min}^{-1})(\|b_0\|_\infty + Lv_{\min}^{-1}\|w_0\|_\infty), \max(0, t - v_{\min}^{-1})\right) \\ + \exp((L+1)(2v_{\min}^{-1} - t))(\|b_0\|_\infty + Lv_{\min}^{-1}\|w_0\|_\infty) \quad (5.88)$$

Using (5.76), (5.84), (5.88) and (5.80), we get for all $\sigma > 0$, $t \geq 0$:

$$\|b[t]\|_\infty + \|w[t]\|_\infty \leq P\left(\exp(Lv_{\min}^{-1})(\|b_0\|_\infty + Lv_{\min}^{-1}\|w_0\|_\infty), \max(0, t - c^{-1} - v_{\min}^{-1})\right) \\ + \exp((L+1)(2v_{\min}^{-1} + c^{-1} - t))(\|b_0\|_\infty + Lv_{\min}^{-1}\|w_0\|_\infty) + \exp(-\sigma(ct-1))\|b_0\|_\infty \\ + \exp(-\sigma(t - v_{\min}^{-1}))\|w_0\|_\infty \quad (5.89)$$

Estimate (5.89) in conjunction with the fact that $L := L(\|w_0\|_\infty + \|b_0\|_\infty)$, where $L: \Re_+ \to \Re_+$ is a non-decreasing function and definition (5.79) of $v_{\min}$ implies that there exists $G \in KL$ such that the following estimate holds for all $t \geq 0$:

$$\|b[t]\|_\infty + \|w[t]\|_\infty \leq G(\|b_0\|_\infty + \|w_0\|_\infty, t) \quad (5.90)$$

Estimate (3.4) for appropriate $Q \in KL$ is a consequence of (5.90) and (2.7), which gives (in conjunction with definition (5.78))

$$\ln\left(\frac{\rho(t,x)}{\rho_{eq}}\right) = w(t,x) + \ln\left(\frac{c+f(\rho_{eq})}{c+v(t,x)}\right) = w(t,x) + \ln\left(1 + f(\rho_{eq})\frac{1-\exp(b(t,x))}{c+f(\rho_{eq})\exp(b(t,x))}\right)$$

and implies the inequalities for all $(t,x) \in \Re_+ \times [0,1]$

$$\left|\ln\left(\frac{\rho(t,x)}{\rho_{eq}}\right)\right| \leq |w(t,x)| + f(\rho_{eq})\left|\frac{1-\exp(b(t,x))}{c+f(\rho_{eq})\exp(b(t,x))}\right| \leq |w(t,x)| + c^{-1}f(\rho_{eq})|1-\exp(b(t,x))|$$
$$\leq |w(t,x)| + c^{-1}f(\rho_{eq})\exp(|b(t,x)|)|b(t,x)|$$

$$|w(t,x)| \leq \left|\ln\left(\frac{\rho(t,x)}{\rho_{eq}}\right)\right| + f(\rho_{eq})\left|\frac{1-\exp(b(t,x))}{c+f(\rho_{eq})\exp(b(t,x))}\right| \leq \left|\ln\left(\frac{\rho(t,x)}{\rho_{eq}}\right)\right| + c^{-1}f(\rho_{eq})|1-\exp(b(t,x))|$$
$$\leq \left|\ln\left(\frac{\rho(t,x)}{\rho_{eq}}\right)\right| + c^{-1}f(\rho_{eq})\exp(|b(t,x)|)|b(t,x)|$$

Indeed, the two above inequalities imply the existence of a function $\varphi \in K_\infty$ such that the estimates

$$\max_{0 \leq x \leq 1}\left(\left|\ln\left(\frac{\rho(t,x)}{\rho_{eq}}\right)\right|\right) \leq \varphi(\|b[t]\|_\infty + \|w[t]\|_\infty), \quad \|w[t]\|_\infty \leq \varphi\left(\|b[t]\|_\infty + \max_{0 \leq x \leq 1}\left(\left|\ln\left(\frac{\rho(t,x)}{\rho_{eq}}\right)\right|\right)\right)$$ hold for all $t \geq 0$.

The proof is complete.   ◁

## 6. Concluding Remarks

The paper provides results for a non-standard, hyperbolic traffic flow model on a bounded domain. The model has been developed for relatively crowded roads and consists of two first-order, hyperbolic PDEs with a dynamic boundary condition, which involves the time derivative of the velocity. Although simple, the proposed model has features that are important from a traffic-theoretic point of view: it is completely anisotropic, i.e., the velocity depends only on the velocity of downstream vehicles, and is a hyperbolic model for which information travels forward exactly at the same speed as traffic. It has been shown that for all physically meaningful initial conditions the model admits a globally defined, unique, classical solution that remains positive and bounded for all times (Theorem 2.1 and Remark 2.2). Moreover, it has been shown that global stabilization in the sup-norm of the logarithmic deviation of the state from its equilibrium point can be achieved for arbitrary equilibria by means of an explicit boundary feedback law which adjusts continuously the inlet flow (Theorem 3.1). It is important to notice that the stabilizing feedback law depends *only* on the inlet velocity. Therefore the measurement requirements for the implementation of the proposed boundary feedback law are minimal. The efficiency of the proposed boundary feedback was demonstrated by means of a numerical example.

Future work may involve the development of more complicated models, retaining the important characteristics of the proposed model, to capture secondary features of traffic flow dynamics. Another direction for future research is the use of sampled-data boundary feedback boundary for the stabilization of unstable equilibria.

## Acknowledgments

Iasson Karafyllis and Markos Papageorgiou were supported by the funding from the European Research Council under the European Union's Seventh Framework Programme (FP/2007-2013) / ERC Grant Agreement n. [321132], project TRAMAN21. Nikolaos Bekiaris-Liberis was supported by the funding from the European Union's Horizon 2020 research and innovation programme under the Marie Sklodowska-Curie grant agreement No 747898, project PADECOT.

# Appendix

**Proof of Lemma 5.1:** Estimate (5.2) is proved by induction for $a=0$ (using (5.1)). For $a>0$ we prove by induction (using (5.1)) the following estimate:

$$x(k) \leq (1+a)^k \left( x(0) + \frac{p}{a+c} \right) + b \frac{(1+a)^k - 1}{a}, \text{ for all } k = 0,1,...,m \qquad (A.1)$$

Using the fact that $\exp(a) \geq 1+a$, we obtain from (A.1) the following estimate:

$$x(k) \leq \exp(ka) \left( x(0) + \frac{p}{a+c} \right) + b \frac{\exp(ka) - 1}{a}, \text{ for all } k = 0,1,...,m \qquad (A.2)$$

Estimate (A.2) in conjunction with the fact that $\exp(ka) - 1 \leq ak \exp(ka)$ implies estimate (5.2). The proof is complete. ◁

**Proof of Proposition 5.2:** Let $T > 0$ be given (arbitrary). We construct a solution $w \in C^1(\Re_+ \times [0,1])$ of the initial-boundary value problem (5.3), (5.4), (5.5), which has Lipschitz derivatives on $[0,T] \times [0,1]$. We follow the methodology of finite-differences presented in the book [14].

Since $v(t,0) > 0$ for all $t \geq 0$, by continuity of $v$ there exists $\varepsilon(T) > 0$ such that $v(t,x) > 0$ for all $t \in [0,T]$, $x \in [0, \varepsilon(T)]$. Let $N^* > 1$ be an integer for which the inequalities

$$2 \max\left\{ \left|\frac{\partial v}{\partial t}(t,x)\right| : (t,x) \in [0,T] \times [0,1] \right\} \leq N^* v_{\min}^2 \text{ and } 2 \leq N^* \varepsilon(T) \tag{A.3}$$

hold, where $v_{\min} := \min\{v(t,x) : (t,x) \in [0,T] \times [0, \varepsilon(T)]\}$. Let $N \geq N^*$ be an integer and consider the parameterized discrete-time system

$$w_i((k+1)\delta) = (1 - \lambda v_i(k\delta))w_i(k\delta) + \lambda v_i(k\delta) w_{i-1}(k\delta), \text{ for } i = 1,\ldots,N, \; k = 0,1,\ldots,m-1 \tag{A.4}$$

$$w_0(k\delta) = a(k\delta), \text{ for } k = 0,1,\ldots,m \tag{A.5}$$

$$w_i(0) = \varphi(ih), \text{ for } i = 1,\ldots,N \tag{A.6}$$

where

$$h := 1/N, \; \delta := \lambda h \tag{A.7}$$

$$v_i(k\delta) := v(k\delta, ih), \text{ for } i = 0,1,\ldots,N, \; k = 0,1,\ldots,m \tag{A.8}$$

$$\lambda := \frac{T}{[Tv_{\max}] + 1} \tag{A.9}$$

$$m := N([Tv_{\max}] + 1) \tag{A.10}$$

$$v_{\max} := \max\{v(t,x) : (t,x) \in [0,T] \times [0,1]\}. \tag{A.11}$$

Notice that the above definitions guarantee that

$$T = m\delta, \tag{A.12}$$

$$\lambda v_{\max} \leq 1. \tag{A.13}$$

Using (A.4), (A.5), (A.8), (A.11) in conjunction with (A.13), we obtain the estimate

$$\max_{i=0,\ldots,N}(|w_i((k+1)\delta)|) \leq \max\left(\|a\|, \max_{i=0,\ldots,N}(|w_i(k\delta)|)\right), \text{ for } k = 0,1,\ldots,m-1 \tag{A.14}$$

where $\|a\| := \max_{0 \leq s \leq t}(|a(s)|)$. It follows from (A.5), (A.6) and (A.14) that the following estimate holds

$$\max_{i=0,\ldots,N}(|w_i(k\delta)|) \leq \max(\|a\|, \|\varphi\|_\infty), \text{ for } k = 0,1,\ldots,m. \tag{A.15}$$

We define the function $w(t,x;N)$ for $(t,x) \in [0,T] \times [0,1]$ and for every integer $N \geq N^*$ (recall that $h = N^{-1}$, $\delta = \lambda h$, $m\delta = T$):

$$w(k\delta, x; N) = (i+1-xN)w_i(k\delta) + (xN-i)w_{i+1}(k\delta)$$
$$\text{with } i = [xN], \text{ for } x \in [0,1), \; k = 0,\ldots,m, \tag{A.16}$$

$$w(k\delta, 1; N) = w_N(k\delta), \text{ for } k = 0,\ldots,m, \tag{A.17}$$

$$w(t,x;N) = (k+1-\lambda^{-1}tN)w(k\delta, x; N) + (\lambda^{-1}tN - k)w((k+1)\delta, x; N)$$
$$\text{with } k = [\lambda^{-1}tN] \text{ for } x \in [0,1], \; t \in [0,T). \tag{A.18}$$

It follows from (A.15) and definitions (A.16), (A.17), (A.18) that the following estimate holds for every integer $N \geq N^*$:

$$\max_{0 \leq x \leq 1}(|w(t,x;N)|) \leq \max(\|a\|, \|\varphi\|_\infty), \text{ for } t \in [0,T]. \tag{A.19}$$

Next define

$$y_i(k\delta) := h^{-1}(w_{i+1}(k\delta) - w_i(k\delta)), \text{ for } i = 0,1,\ldots,N-1, \ k = 0,1,\ldots,m, \quad (A.20)$$

$$y_N(k\delta) := y_{N-1}(k\delta), \text{ for } k = 0,1,\ldots,m. \quad (A.21)$$

Equations (A.4), (A.5) in conjunction with definitions (A.20), (A.21) imply that the following equalities hold:

$$y_i((k+1)\delta) = (1 - \lambda v_{i+1}(k\delta))y_i(k\delta) + \lambda v_{i+1}(k\delta)y_{i-1}(k\delta) - \lambda(v_{i+1}(k\delta) - v_i(k\delta))y_{i-1}(k\delta),$$
$$\text{for } i = 1,\ldots,N-1, \ k = 0,1,\ldots,m-1, \quad (A.22)$$

$$y_0((k+1)\delta) = (1 - \lambda v_1(k\delta))y_0(k\delta) - \lambda \frac{a((k+1)\delta) - a(k\delta)}{\delta}, \text{ for } k = 0,1,\ldots,m-1 \quad (A.23)$$

Using the fact that $|v_{i+1}(k\delta) - v_i(k\delta)| \leq h\|v_x\|$ for all $i = 0,\ldots,N-1$, $k = 0,1,\ldots,m$, where $\|v_x\| := \max\left\{\left|\frac{\partial v}{\partial x}(t,x)\right| : (t,x) \in [0,T] \times [0,1]\right\}$ (a direct consequence of definition (A.8)), the fact $|a((k+1)\delta) - a(k\delta)| \leq \delta\|\dot a\|$ for all $k = 0,1,\ldots,m-1$, where $\|\dot a\| := \max\{|\dot a(t)| : t \in [0,T]\}$ and the fact that $v_1(k\delta) \geq v_{\min} > 0$ for all $k = 0,1,\ldots,m$, where $v_{\min} := \min\{v(t,x) : (t,x) \in [0,T] \times [0,\varepsilon(T)]\}$ (a consequence of (A.3), (A.7), (A.8) which imply that $2h \leq \varepsilon(T)$), in conjunction with (A.21), (A.22), (A.23), (A.7), (A.8), (A.11), (A.13), we get for $k = 0,1,\ldots,m-1$:

$$\max_{i=0,\ldots,N}(|y_i((k+1)\delta)|) \leq \max\left((1 + \delta\|v_x\|)\max_{i=0,\ldots,N}(|y_i(k\delta)|), (1 - \lambda v_{\min})\max_{i=0,\ldots,N}(|y_i(k\delta)|) + \lambda\|\dot a\|\right) \quad (A.24)$$

Using (A.24) in conjunction with (A.12), the fact that $|y_i(0)| \leq \|\varphi'\|_\infty$ for $i = 0,\ldots,N$ (a direct consequence of definitions (A.5), (A.6), (A.20), (A.21) and the fact that $a(0) = \varphi(0)$) and Lemma 5.1, we obtain the estimate:

$$\max_{i=0,\ldots,N}(|y_i(k\delta)|) \leq Y := \exp(T\|v_x\|)\left(\|\varphi'\|_\infty + \frac{\|\dot a\|}{v_{\min}}\right), \text{ for } k = 0,1,\ldots,m \quad (A.25)$$

Next define

$$p_i(k\delta) := \delta^{-1}(w_i((k+1)\delta) - w_i(k\delta)), \text{ for } i = 0,1,\ldots,N, \ k = 0,1,\ldots,m-1, \quad (A.26)$$

$$p_i(m\delta) := p_i((m-1)\delta), \text{ for } i = 0,1,\ldots,N. \quad (A.27)$$

Using (A.4), (A.5), (A.7), (A.11), (A.20), (A.26), (A.27) (which imply that $p_i(k\delta) = -v_i(k\delta)y_{i-1}(k\delta)$, for $i = 1,\ldots,N$, $k = 0,1,\ldots,m-1$ and $p_0(k\delta) = \frac{a((k+1)\delta) - a(k\delta)}{\delta}$, for $k = 0,1,\ldots,m-1$), the fact $|a((k+1)\delta) - a(k\delta)| \leq \delta\|\dot a\|$ for all $k = 0,1,\ldots,m-1$, where $\|\dot a\| := \max\{|\dot a(t)| : t \in [0,T]\}$ and estimate (A.25), we obtain

$$\max_{i=0,\ldots,N}(|p_i(k\delta)|) \leq \max(v_{\max}Y, \|\dot a\|), \text{ for } k = 0,1,\ldots,m \quad (A.28)$$

Definitions (A.20), (A.26) in conjunction with estimates (A.25), (A.28) imply the estimate for $i,j = 0,1,\ldots,N$, $k,l = 0,1,\ldots,m$:

$$|w_i(k\delta) - w_j(l\delta)| \leq h|i-j|Y + \delta|k-l|\max(v_{\max}Y, \|\dot a\|) \quad (A.29)$$

Estimate (A.29) in conjunction with definitions (A.16), (A.17), (A.18) implies that there exists a constant $L_1 = L_1(T, a, v, \varphi) > 0$ such that the following Lipschitz inequality holds for every integer $N \geq N^*$:

$$|w(t,x;N) - w(\tau,z;N)| \leq L_1(|x-z| + |t-\tau|), \text{ for all } t,\tau \in [0,T] \text{ and } x,z \in [0,1]. \quad (A.30)$$

Next define:

$$\psi(k\delta) := h^{-1}\left(\frac{a((k+1)\delta) - a(k\delta)}{\delta} + v_1(k\delta)y_0(k\delta)\right), \text{ for } k = 0,1,\ldots,m-1 \quad (A.31)$$

It follows from (A.31) and (A.23) that the following equation holds for $k = 0,1,\ldots,m-2$:

$$\psi((k+1)\delta) = \frac{v_1((k+1)\delta)}{v_1(k\delta)}(1 - \lambda v_1(k\delta))\psi(k\delta) + \lambda \frac{a((k+2)\delta) - 2a((k+1)\delta) + a(k\delta)}{\delta^2}$$
$$- \lambda \frac{1}{v_1(k\delta)} \frac{a((k+1)\delta) - a(k\delta)}{\delta} \frac{v_1((k+1)\delta) - v_1(k\delta)}{\delta} \quad (A.32)$$

It should be noticed that inequalities (A.3) in conjunction with (A.7), (A.8) and definition $v_{\min} := \min\{v(t,x) : (t,x) \in [0,T] \times [0, \varepsilon(T)]\}$ implies the following inequality for $k = 0,1,\ldots,m-1$:

$$\frac{v_1((k+1)\delta)}{v_1(k\delta)}(1 - \lambda v_1(k\delta)) \leq 1 - \lambda \frac{v_{\min}}{2} \quad (A.33)$$

It follows from (A.32), (A.33) in conjunction with the fact that $|a((k+1)\delta) - a(k\delta)| \leq \delta \|\dot{a}\|$ for all $k = 0,1,\ldots,m-1$, where $\|\dot{a}\| := \max\{|\dot{a}(t)| : t \in [0,T]\}$, $|v_1((k+1)\delta) - v_1(k\delta)| \leq \delta \|v_t\|$ for all $k = 0,1,\ldots,m-1$, where $\|v_t\| := \max\left\{\left|\frac{\partial v}{\partial t}(t,x)\right| : (t,x) \in [0,T] \times [0,1]\right\}$ (a direct consequence of definition (A.8)), definition $v_{\min} := \min\{v(t,x) : (t,x) \in [0,T] \times [0, \varepsilon(T)]\}$, the fact that $v_1(k\delta) \geq v_{\min} > 0$ for all $k = 0,1,\ldots,m$ (a consequence of (A.3), (A.7), (A.8) which imply that $2h \leq \varepsilon(T)$), the fact that $|a((k+2)\delta) - 2a((k+1)\delta) + a(k\delta)| \leq \delta^2 \|\ddot{a}\|$ for all $k = 0,1,\ldots,m-2$, where $\|\ddot{a}\| := ess \sup\{|\ddot{a}(t)| : t \in [0,T]\}$, that the following inequality holds for $k = 0,1,\ldots,m-2$:

$$|\psi((k+1)\delta)| \leq \left(1 - \lambda \frac{v_{\min}}{2}\right)|\psi(k\delta)| + \lambda \|\ddot{a}\| + \lambda \frac{\|\dot{a}\| \|v_t\|}{v_{\min}} \quad (A.34)$$

Consequently, we obtain directly (by induction) the following estimate for $k = 0,1,\ldots,m-1$:

$$|\psi(k\delta)| \leq |\psi(0)| + 2\frac{\|\dot{a}\| \|v_t\| + v_{\min} \|\ddot{a}\|}{v_{\min}^2} \quad (A.35)$$

Definitions (A.5), (A.6), (A.8), (A.11), (A.20), (A.31) in conjunction with the facts that $a(0) = \varphi(0)$, $\dot{a}(0) + v(0,0)\varphi'(0) = 0$, definitions $\|v_x\| := \max\left\{\left|\frac{\partial v}{\partial x}(t,x)\right| : (t,x) \in [0,T] \times [0,1]\right\}$, $\|\ddot{a}\| := ess \sup\{|\ddot{a}(t)| : t \in [0,T]\}$, $\|\varphi''\|_\infty := ess \sup\{|\varphi''(x)| : x \in [0,1]\}$ and inequality (A.13) implies that:

$$|\psi(0)| \leq \frac{\|\ddot{a}\|}{2v_{\max}} + |\varphi'(0)| \|v_x\| + v_{\max} \frac{\|\varphi''\|_\infty}{2} \quad (A.36)$$

Consequently, we get from (A.35), (A.36) for $k = 0,1,\ldots,m-1$:

$$|\psi(k\delta)| \leq M := \frac{\|\ddot{a}\|}{2v_{\max}} + |\varphi'(0)| \|v_x\| + v_{\max} \frac{\|\varphi''\|_\infty}{2} + 2\frac{\|\dot{a}\| \|v_t\| + v_{\min} \|\ddot{a}\|}{v_{\min}^2} \quad (A.37)$$

We define the function $y(t,x;N)$ for $(t,x) \in [0,T] \times [0,1]$ and for every integer $N \geq N^*$ (recall that $h = N^{-1}$, $\delta = \lambda h$, $m\delta = T$):

$$y(k\delta, x; N) = (i + 1 - xN)y_i(k\delta) + (xN - i)y_{i+1}(k\delta)$$
$$\text{with } i = [xN], \text{ for } x \in [0,1), \; k = 0,\ldots,m, \quad (A.38)$$

$$y(k\delta, 1; N) = y_N(k\delta), \text{ for } k = 0,\ldots,m, \quad (A.39)$$

$$y(t,x;N) = (k + 1 - \lambda^{-1}tN)y(k\delta, x; N) + (\lambda^{-1}tN - k)y((k+1)\delta, x; N)$$
$$\text{with } k = [\lambda^{-1}tN] \text{ for } x \in [0,1], \; t \in [0,T). \quad (A.40)$$

It follows from (A.25) and definitions (A.38), (A.39), (A.40) that the following estimate holds for every integer $N \geq N^*$:

$$\max_{0 \leq x \leq 1}(|y(t,x;N)|) \leq Y := \exp(T\|v_x\|)\left(\|\varphi'\|_\infty + \frac{\|\dot{a}\|}{v_{\min}}\right), \text{ for } t \in [0,T]. \quad (A.41)$$

We also define the function $p(t, x; N)$ for $(t, x) \in [0, T] \times [0, 1]$ and for every integer $N \geq N^*$ (recall that $h = N^{-1}$, $\delta = \lambda h$, $m\delta = T$):

$$p(k\delta, x; N) = (i + 1 - xN) p_i(k\delta) + (xN - i) p_{i+1}(k\delta)$$
with $i = [xN]$, for $x \in [0, 1)$, $k = 0, \ldots, m$, \quad (A.42)

$$p(k\delta, 1; N) = p_N(k\delta), \text{ for } k = 0, \ldots, m,$$ \quad (A.43)

$$p(t, x; N) = (k + 1 - \lambda^{-1} tN) p(k\delta, x; N) + (\lambda^{-1} tN - k) p((k+1)\delta, x; N)$$
with $k = [\lambda^{-1} tN]$ for $x \in [0, 1]$, $t \in [0, T)$. \quad (A.44)

It follows from (A.28) and definitions (A.42), (A.43), (A.44) that the following estimate holds for every integer $N \geq N^*$:

$$\max_{0 \leq x \leq 1} (|p(t, x; N)|) \leq \max(v_{\max} Y, \|\dot{a}\|), \text{ for } t \in [0, T].$$ \quad (A.45)

Next define

$$\omega_i(k\delta) := h^{-1}(y_{i+1}(k\delta) - y_i(k\delta)), \text{ for } i = 0, 1, \ldots, N-1, \, k = 0, 1, \ldots, m.$$ \quad (A.46)

Definitions (A.21), (A.31), (A.46) in conjunction with equations (A.22), (A.23) implies that the following equalities hold:

$$\omega_i((k+1)\delta) = (1 - \lambda v_{i+2}(k\delta)) \omega_i(k\delta) + \lambda v_{i+2}(k\delta) \omega_{i-1}(k\delta)$$
$$- 2\lambda (v_{i+2}(k\delta) - v_{i+1}(k\delta)) \omega_{i-1}(k\delta) - \lambda \frac{v_{i+2}(k\delta) - 2v_{i+1}(k\delta) + v_i(k\delta)}{h} y_{i-1}(k\delta),$$
for $i = 1, \ldots, N-2$, $k = 0, 1, \ldots, m-1$ \quad (A.47)

$$\omega_{N-1}(k\delta) = 0, \text{ for } k = 0, 1, \ldots, m.$$ \quad (A.48)

$$\omega_0((k+1)\delta) = (1 - \lambda v_2(k\delta)) \omega_0(k\delta) - \lambda \frac{v_2(k\delta) - v_1(k\delta)}{h} y_0(k\delta) + \lambda \psi(k\delta), \text{ for } k = 0, 1, \ldots, m-1$$ \quad (A.49)

Using the facts that $|v_{i+2}(k\delta) - 2v_{i+1}(k\delta) + v_i(k\delta)| \leq h^2 \|v_{xx}\|$, $|v_{i+2}(k\delta) - v_{i+1}(k\delta)| \leq h \|v_x\|$ for $i = 0, \ldots, N-2$, $k = 0, 1, \ldots, m$,

where $\|v_{xx}\| := \sup \left\{ \frac{\left| \frac{\partial v}{\partial x}(t, x) - \frac{\partial v}{\partial x}(t, z) \right|}{|x - z|} : (t, x, z) \in [0, T] \times [0, 1]^2, x \neq z \right\}$, $\|v_x\| := \max \left\{ \left| \frac{\partial v}{\partial x}(t, x) \right| : (t, x) \in [0, T] \times [0, 1] \right\}$

(a direct consequence of definition (A.8)) and the fact that $v_2(k\delta) \geq v_{\min} > 0$ for all $k = 0, 1, \ldots, m$, where $v_{\min} := \min \{v(t, x) : (t, x) \in [0, T] \times [0, \varepsilon(T)]\}$ (a consequence of (A.3), (A.7), (A.8) which imply that $2h \leq \varepsilon(T)$), in conjunction with (A.25), (A.37), (A.8), (A.11), (A.13), we get for $k = 0, 1, \ldots, m-1$:

$$\max_{i=0, \ldots N-1} (|\omega_i((k+1)\delta)|) \leq \max \left( (1 + 2\delta \|v_x\|) \max_{i=0, \ldots N-1} (|\omega_i(k\delta)|) + \delta \|v_{xx}\| Y, (1 - \lambda v_{\min}) \max_{i=0, \ldots N-1} (|\omega_i(k\delta)|) + \lambda (\|v_x\| Y + M) \right)$$ \quad (A.50)

It follows from (A.12) and Lemma 5.1 for $k = 0, 1, \ldots, m$:

$$\max_{i=0, \ldots N-1} (|\omega_i((k\delta)|) \leq \exp(2T \|v_x\|) \left( \max_{i=0, \ldots N-1} (|\omega_i(0)|) + \frac{\|v_x\| Y + M}{v_{\min}} + T \|v_{xx}\| Y \right)$$ \quad (A.51)

Definition (A.5), (A.6), (A.46), (A.20), (A.21) in conjunction with the fact that $a(0) = \varphi(0)$ implies that $|\omega_i(0)| \leq \|\varphi''\|_\infty$ for all $i = 0, \ldots, N-1$, where $\|\varphi''\|_\infty := \text{ess sup}\{|\varphi''(x)| : x \in [0, 1]\}$. Therefore, it follows from (A.51) that the following estimate holds for $k = 0, 1, \ldots, m$:

$$\max_{i=0, \ldots N-1} (|\omega_i((k\delta)|) \leq \Omega := \exp(2T \|v_x\|) \left( \|\varphi''\|_\infty + \frac{\|v_x\| Y + M}{v_{\min}} + T \|v_{xx}\| Y \right)$$ \quad (A.52)

Next define

$$\eta_i(k\delta) := \delta^{-1}(y_i((k+1)\delta) - y_i(k\delta)), \text{ for } i = 0, 1, \ldots, N, \, k = 0, 1, \ldots, m-1.$$ \quad (A.53)

Equations (A.22), (A.23) in conjunction with definitions (A.31), (A.46), (A.53), imply the relations

$$\eta_i(k\delta) = -v_{i+1}(k\delta)\omega_{i-1}(k\delta) - \frac{v_{i+1}(k\delta) - v_i(k\delta)}{h} y_{i-1}(k\delta), \text{ for } i = 1,\ldots,N-1,\ k = 0,1,\ldots,m-1 \quad (A.54)$$

$$\eta_0(k\delta) = -\psi(k\delta), \text{ for } k = 0,1,\ldots,m-1 \quad (A.55)$$

which combined with (A.11), (A.21), (A.52), (A.37) and the fact that $|v_{i+1}(k\delta) - v_i(k\delta)| \leq h\|v_x\|$ for all $i = 0,\ldots,N-1$, $k = 0,1,\ldots,m$, where, $\|v_x\| := \max\left\{\left|\frac{\partial v}{\partial x}(t,x)\right| : (t,x) \in [0,T] \times [0,1]\right\}$ (a direct consequence of definition (A.8)) give:

$$\max_{i=0,\ldots,N}\left(|\eta_i((k\delta))|\right) \leq v_{\max}\Omega + \|v_x\|Y + M, \text{ for } k = 0,1,\ldots,m-1 \quad (A.56)$$

Definitions (A.46), (A.53) in conjunction with estimates (A.52), (A.56) imply the following estimate for $i,j = 0,1,\ldots,N$, $k,l = 0,1,\ldots,m$:

$$|y_i(k\delta) - y_j(l\delta)| \leq h|i-j|\Omega + \delta|k-l|(v_{\max}\Omega + M + Y\|v_x\|) \quad (A.57)$$

Estimate (A.57) in conjunction with definitions (A.38), (A.39), (A.40) implies that there exists a constant $L_2 = L_2(T,a,v,\varphi) > 0$ such that the following Lipschitz inequality holds for every integer $N \geq N^*$:

$$|y(t,x;N) - y(\tau,z;N)| \leq L_2(|x-z| + |t-\tau|), \text{ for all } t,\tau \in [0,T] \text{ and } x,z \in [0,1]. \quad (A.58)$$

Next define

$$\zeta_i(k\delta) := h^{-1}(p_{i+1}(k\delta) - p_i(k\delta)), \text{ for } i = 0,1,\ldots,N-1,\ k = 0,1,\ldots,m. \quad (A.59)$$

Since $p_i(k\delta) = -v_i(k\delta)y_{i-1}(k\delta)$, for $i = 1,\ldots,N$, $k = 0,1,\ldots,m-1$ and $p_0(k\delta) = \frac{a((k+1)\delta) - a(k\delta)}{\delta}$, for $k = 0,1,\ldots,m-1$, we obtain from (A.27), (A.31), (A.46):

$$\zeta_i(k\delta) = -\frac{v_{i+1}(k\delta) - v_i(k\delta)}{h} y_i(k\delta) - v_i(k\delta)\omega_{i-1}(k\delta), \text{ for } i = 1,\ldots,N-1,\ k = 0,1,\ldots,m-1 \quad (A.60)$$

$$\zeta_0(k\delta) = -\psi(k\delta), \text{ for } k = 0,1,\ldots,m-1 \quad (A.61)$$

$$\zeta_i(m\delta) = \zeta_i((m-1)\delta), \text{ for } i = 0,1,\ldots,N-1 \quad (A.62)$$

Using the facts that $|v_{i+1}(k\delta) - v_i(k\delta)| \leq h\|v_x\|$ for all $i = 0,\ldots,N-1$, $k = 0,1,\ldots,m$, where $\|v_x\| := \max\left\{\left|\frac{\partial v}{\partial x}(t,x)\right| : (t,x) \in [0,T] \times [0,1]\right\}$ (a direct consequence of definition (A.8)) in conjunction with (A.60), (A.61), (A.62), (A.11), (A.25), (A.37), (A.52), we get the following estimate for $k = 0,1,\ldots,m$:

$$\max_{i=0,\ldots,N-1}\left(|\zeta_i((k\delta))|\right) \leq Z := \|v_x\|Y + v_{\max}\Omega + M \quad (A.63)$$

Finally, define

$$\xi_i(k\delta) := \delta^{-1}(p_i((k+1)\delta) - p_i(k\delta)), \text{ for } i = 0,1,\ldots,N,\ k = 0,1,\ldots,m-1. \quad (A.64)$$

Similarly, as above and using (A.53), we obtain:

$$\xi_i(k\delta) = -\frac{v_i((k+1)\delta) - v_i(k\delta)}{\delta} y_{i-1}((k+1)\delta) - v_i(k\delta)\eta_{i-1}(k\delta), \text{ for } i = 1,\ldots,N,\ k = 0,1,\ldots,m-2 \quad (A.65)$$

$$\xi_0(k\delta) = \frac{a((k+2)\delta) - 2a((k+1)\delta) + a(k\delta)}{\delta^2}, \text{ for } k = 0,1,\ldots,m-2 \quad (A.66)$$

$$\xi_i((m-1)\delta) = 0, \text{ for } i = 0,1,\ldots,N \quad (A.67)$$

It follows from (A.11), (A.25), (A.56), (A.65), (A.66), (A.67) in conjunction with the fact that $|v_i((k+1)\delta) - v_i(k\delta)| \leq \delta \|v_t\|$ for all $i = 0,1,..., N$, $k = 0,1,...,m-1$, where $\|v_t\| := \max\left\{\left|\frac{\partial v}{\partial t}(t,x)\right| : (t,x) \in [0,T] \times [0,1]\right\}$ (a direct consequence of definition (A.8)), the fact that $|a((k+2)\delta) - 2a((k+1)\delta) + a(k\delta)| \leq \delta^2 \|\ddot{a}\|$ for all $k = 0,1,...,m-2$, where $\|\ddot{a}\| := ess\sup\{|\ddot{a}(t)| : t \in [0,T]\}$, that the following inequality holds for $k = 0,1,...,m-1$:

$$\max_{i=0,...N} (|\xi_i((k\delta)|) \leq \Xi := \|v_t\| Y + v_{\max} (\|v_x\| Y + v_{\max} \Omega + M) + \|\ddot{a}\| \tag{A.68}$$

Definitions (A.59), (A.64) in conjunction with estimates (A.63), (A.68) imply the following estimate for $i, j = 0,1,...,N$, $k, l = 0,1,...,m$:

$$|p_i(k\delta) - p_j(l\delta)| \leq h|i-j|Z + \delta|k-l|\Xi \tag{A.69}$$

Estimate (A.69) in conjunction with definitions (A.42), (A.43), (A.44) implies that there exists a constant $L_3 = L_3(T, a, v, \varphi) > 0$ such that the following Lipschitz inequality holds for every integer $N \geq N^*$:

$$|p(t, x; N) - p(\tau, z; N)| \leq L_3(|x-z| + |t-\tau|), \text{ for all } t, \tau \in [0,T] \text{ and } x, z \in [0,1]. \tag{A.70}$$

It follows from (A.19), (A.30), (A.41), (A.58), (A.45), (A.70) that the sequences of functions $\{w(\cdot; N)\}_{N=N^*}^{\infty}$, $\{y(\cdot; N)\}_{N=N^*}^{\infty}$, $\{p(\cdot; N)\}_{N=N^*}^{\infty}$ are uniformly bounded and equicontinuous. Therefore, compactness of $[0,T] \times [0,1]$ and the Arzela-Ascoli theorem implies that there exist Lipschitz functions $w : [0,T] \times [0,1] \to \Re$, $y : [0,T] \times [0,1] \to \Re$, $p : [0,T] \times [0,1] \to \Re$ and subsequences $\{w(\cdot; q)\}_{q=N^*}^{\infty} \subset \{w(\cdot; N)\}_{N=N^*}^{\infty}$, $\{y(\cdot; q)\}_{q=N^*}^{\infty} \subset \{y(\cdot; N)\}_{N=N^*}^{\infty}$, $\{p(\cdot; q)\}_{q=N^*}^{\infty} \subset \{p(\cdot; N)\}_{N=N^*}^{\infty}$ which converge uniformly on $[0,T] \times [0,1]$ to $w, y$ and $p$, respectively. Moreover, the functions $w, y$ and $p$ satisfy the same bounds with $w(\cdot; N), y(\cdot; N)$ and $p(\cdot; N)$, i.e., $\|w[t]\|_{\infty} \leq \max(\|a\|_{\infty}, \|\varphi\|_{\infty})$, $\|y[t]\|_{\infty} \leq \exp(T\|v_x\|)\left(\|\varphi'\|_{\infty} + \frac{\|\dot{a}\|}{v_{\min}}\right)$ and $\|p[t]\|_{\infty} \leq v_{\max} \exp(T\|v_x\|)\left(\|\varphi'\|_{\infty} + \frac{\|\dot{a}\|}{v_{\min}}\right)$, for $t \in [0,T]$.

We remark that in what follows the convergent subsequences $\{w(\cdot; q)\}_{q=N^*}^{\infty} \subset \{w(\cdot; N)\}_{N=N^*}^{\infty}$, $\{y(\cdot; q)\}_{q=N^*}^{\infty} \subset \{y(\cdot; N)\}_{N=N^*}^{\infty}$, $\{p(\cdot; q)\}_{q=N^*}^{\infty} \subset \{p(\cdot; N)\}_{N=N^*}^{\infty}$, will be denoted by $\{w(\cdot; N)\}_{N=N^*}^{\infty}$, $\{y(\cdot; N)\}_{N=N^*}^{\infty}$, $\{p(\cdot; N)\}_{N=N^*}^{\infty}$.

We next show that $p(t, x) + v(t, x) y(t, x) = 0$ for all $(t, x) \in [0,T] \times [0,1]$. Let $(t, x) \in [0,T] \times (0,1]$ be given (arbitrary). Since $p_i(k\delta) = -v_i(k\delta) y_{i-1}(k\delta)$, for $i = 1,...,N$, $k = 0,1,...,m-1$, we obtain using definitions (A.7), (A.38), (A.39), (A.40), (A.42), (A.43), (A.44), inequalities (A.11), (A.25), (A.58), (A.70) and definitions $\|v_t\| := \max\left\{\left|\frac{\partial v}{\partial t}(t,x)\right| : (t,x) \in [0,T] \times [0,1]\right\}$, $\|v_x\| := \max\left\{\left|\frac{\partial v}{\partial x}(t,x)\right| : (t,x) \in [0,T] \times [0,1]\right\}$, for $k = [\lambda^{-1} t N]$ and $i$ being the smallest integer which is greater or equal to $xN$:

$$|p(t, x; N) + v(t, x) y(t, x; N)| \leq |p(t, x; N) - p_i(k\delta)| + v_i(k\delta)|y(t, x; N) - y_{i-1}(k\delta)| + |v(t, x) - v_i(k\delta)| |y(t, x; N)|$$

$$\leq L_3(t - k\delta + ih - x) + L_2 v_{\max}(t - k\delta + x - (i-1)h) + Y\|v_t\|(t - k\delta) + Y\|v_x\|(ih - x)$$

$$\leq (L_3 + L_2 v_{\max} + Y\|v_t\| + Y\|v_x\|)(\delta + h) \leq (L_3 + L_2 v_{\max} + Y\|v_t\| + Y\|v_x\|)(\lambda + 1) N^{-1}$$

In the above derivation, we have used the facts that $t - k\delta \leq \delta$ and $(i-1)h < x \leq ih$. Since $\{y(\cdot; N)\}_{N=N^*}^{\infty}$, $\{p(\cdot; N)\}_{N=N^*}^{\infty}$ converge uniformly to $y$ and $p$ as $N \to +\infty$, the above inequality shows that $p(t, x) + v(t, x) y(t, x) = 0$ for all $(t, x) \in [0,T] \times (0,1]$. Continuity of $p, v, y$ implies that the equality $p(t, x) + v(t, x) y(t, x) = 0$ for all $(t, x) \in [0,T] \times [0,1]$.

We next show that $y(t,x) = \frac{\partial w}{\partial x}(t,x)$ for all $(t,x) \in [0,T] \times [0,1]$. Equivalently, we show that $w(t,x) - w(t,0) = \int_0^x y(t,z) dz$, for all $(t,x) \in [0,T] \times [0,1]$. Let $(t,x) \in [0,T] \times [0,1]$ be given (arbitrary). Using definitions (A.7), (A.16), (A.17), (A.18), (A.38), (A.39), (A.40), inequalities (A.30), (A.58), (A.41), we obtain for $k = [\lambda^{-1} tN]$ and $i = [xN]$:

$$\left| w(t,x;N) - w(t,0;N) - \int_0^x y(t,z;N) dz \right| \leq |w(t,x;N) - w(k\delta, x;N)| + |w(k\delta, x;N) - w(k\delta, ih;N)|$$

$$+ |w(k\delta, 0;N) - w(t,0;N)| + \left| \int_{ih}^x y(t,z;N) dz \right| + \left| w_i(k\delta) - w_0(k\delta) - \int_0^{ih} y(t,z;N) dz \right|$$

$$\leq L_1(2t - 2k\delta + x - ih) + Y(x - ih) + \left| h \sum_{s=0}^{i-1} y_s(k\delta) - \int_0^{ih} y(t,z;N) dz \right|$$

$$\leq L_1(2\lambda + 1)h + Yh + \left| \sum_{s=0}^{i-1} \int_{sh}^{(s+1)h} y(k\delta, sh;N) dz - \int_0^{ih} y(t,z;N) dz \right|$$

$$\leq (L_1(2\lambda + 1) + Y)h + \left| \sum_{s=0}^{i-1} \int_{sh}^{(s+1)h} (y(k\delta, sh;N) - y(t,z;N)) dz \right|$$

$$\leq (L_1(2\lambda + 1) + Y)h + \sum_{s=0}^{i-1} \int_{sh}^{(s+1)h} |y(k\delta, sh;N) - y(t,z;N)| dz$$

$$\leq (L_1(2\lambda + 1) + Y)h + (L_2(t - k\delta) + L_2 h) \sum_{s=0}^{i-1} \int_{sh}^{(s+1)h} dz = (L_1(2\lambda + 1) + Y)h + (L_2(t - k\delta) + L_2 h)ih$$

$$\leq (L_1(2\lambda + 1) + Y + L_2(\lambda + 1))h = (L_1(2\lambda + 1) + Y + L_2(\lambda + 1))N^{-1}$$

In the above derivation, we have used the facts that $t - k\delta \leq \delta$ and $x - ih \leq h$. Since $\{y(\cdot;N)\}_{N=N^*}^\infty$, $\{w(\cdot;N)\}_{N=N^*}^\infty$ converge uniformly to $y$ and $w$ as $N \to +\infty$, the above inequality shows that $w(t,x) - w(t,0) = \int_0^x y(t,z) dz$ for all $(t,x) \in [0,T] \times [0,1]$.

Using exactly the same approach we show that $p(t,x) = \frac{\partial w}{\partial t}(t,x)$ for all $(t,x) \in [0,T] \times [0,1]$ (i.e., equivalently, we show that $w(t,x) - w(0,x) = \int_0^t p(\tau, x) d\tau$, for all $(t,x) \in [0,T] \times [0,1]$).

Uniqueness follows from a contradiction argument. Suppose that there exist two solutions $w, \tilde{w} \in C^1(\Re_+ \times [0,1])$ of the initial-boundary value problem (5.3), (5.4), (5.5). It then follows that the function $e = w - \tilde{w}$ satisfies the following equations:

$$\frac{\partial e}{\partial t}(t,x) + v(t,x) \frac{\partial e}{\partial x}(t,x) = 0, \text{ for } t \geq 0, \ x \in [0,1] \tag{A.71}$$

$$e(0,x) = 0, \text{ for } x \in [0,1] \tag{A.72}$$

$$e(t,0) = 0, \text{ for } t \geq 0 \tag{A.73}$$

Using the functional $V(t) = \int_0^1 e^2(t,x) dx$ on $[0,T]$ for arbitrary $T > 0$, we have by virtue of (A.71) and (A.73) for every $t \in [0,T]$:

$$\dot{V}(t) = 2\int_0^1 e(t,x)\frac{\partial e}{\partial t}(t,x)dx = -2\int_0^1 v(t,x)e(t,x)\frac{\partial e}{\partial x}(t,x)dx = -\int_0^1 v(t,x)\frac{\partial}{\partial x}\left(e^2(t,x)\right)dx$$

$$= -v(t,1)e^2(t,1) + \int_0^1 e^2(t,x)\frac{\partial v}{\partial x}(t,x)dx \leq \|v_x\|V(t)$$

where $\|v_x\| := \max\left\{\left|\frac{\partial v}{\partial x}(t,x)\right| : (t,x) \in [0,T] \times [0,1]\right\}$. Gronwall's lemma implies that $V(t) \leq \exp(\|v_x\|t)V(0)$ for all $t \in [0,T]$ and consequently (using (A.72)), we get $V(t) = 0$ for all $t \in [0,T]$. This implies $w \equiv \tilde{w}$.

Finally, we assume that there exists a constant $v_{\min} > 0$ such that $v(t,x) \geq v_{\min}$ for all $t \geq 0$, $x \in [0,1]$ and that $a \equiv 0$. Let $T > v_{\min}^{-1}$ be given (arbitrary). Consider the parameterized family of functionals $V_\sigma(t) = \int_0^1 \exp(-\sigma x)w^2(t,x)dx$ on $[0,T]$ with parameter $\sigma > 0$. Using (5.3) and (5.5) with $a \equiv 0$ and the fact that $v(t,x) \geq v_{\min}$ for all $t \geq 0$, we get:

$$\dot{V}_\sigma(t) = 2\int_0^1 \exp(-\sigma x)w(t,x)\frac{\partial w}{\partial t}(t,x)dx = -2\int_0^1 \exp(-\sigma x)v(t,x)w(t,x)\frac{\partial w}{\partial x}(t,x)dx$$

$$= -\int_0^1 \exp(-\sigma x)v(t,x)\frac{\partial}{\partial x}\left(w^2(t,x)\right)dx$$

$$= -v(t,1)\exp(-\sigma)w^2(t,1) + \int_0^1 \exp(-\sigma x)w^2(t,x)\frac{\partial v}{\partial x}(t,x)dx - \sigma\int_0^1 \exp(-\sigma x)w^2(t,x)v(t,x)dx$$

$$\leq -(\sigma v_{\min} - \|v_x\|)V_\sigma(t)$$

It follows that $V_\sigma(t) \leq \exp(-(\sigma v_{\min} - \|v_x\|)t)V_\sigma(0)$, for all $t \in [0,T]$. Since $\exp(-\sigma)\|w[t]\|_2^2 \leq V_\sigma(t) \leq \|w[t]\|_2^2$ for all $t \in [0,T]$, the previous inequality implies the estimate $\|w[t]\|_2^2 \leq \exp\left(-\sigma\left((v_{\min} - \sigma^{-1}\|v_x\|)t - 1\right)\right)\|\varphi\|_2^2$, for all $t \in [0,T]$. Since $\lim_{\sigma \to +\infty}\left(-\sigma\left((v_{\min} - \sigma^{-1}\|v_x\|)t - 1\right)\right) = -\infty$ for each fixed $t \in (v_{\min}^{-1}, T]$, we conclude that $\|w[t]\|_2 = 0$, for all $t \in (v_{\min}^{-1}, T]$. Therefore, by continuity of $w$ and since $T > v_{\min}^{-1}$ is arbitrary, we conclude that $w(t,x) = 0$ for all $x \in [0,1]$ and $t \geq v_{\min}^{-1}$. The proof is complete. ◁